\documentclass[reqno,draft,intlimits]{amsart}

\usepackage[active]{srcltx}

\usepackage[reqno]{amsmath}
\usepackage{amssymb}

\usepackage[arrow,matrix]{xy}

\usepackage{wrapfig}

\makeatletter
\long\def\@makecaption#1#2{%
  \vskip3pt
  \sbox\@tempboxa{\small#1. #2}%
    \ifdim \wd\@tempboxa >\hsize
    \small#1. #2\par
  \else
    \global \@minipagefalse
    \hb@xt@\hsize{\hfil\box\@tempboxa\hfil}%
  \fi
  \vskip0pt}
\makeatother

\newtheorem{thm}{Theorem}[section] 
  
\newtheorem{proc}{Proposition}[section]

\newtheorem{lem}{Lemma}[section]

\newtheorem{defi}{Definition}[section]
\newtheorem{rmk}{Remark}[section]
\newtheorem{ex} {Example}[section]

\input cyracc.def
\font \tencyr=wncyr10
\newfam\cyrfam
\font\tencyr=wncyr10
\font\sevencyr=wncyr7
\font\fivecyr=wncyr5
\def\cyr{\fam\cyrfam\tencyr\cyracc}

\textfont\cyrfam=\tencyr \scriptfont\cyrfam=\sevencyr
\scriptscriptfont\cyrfam=\fivecyr

\newcommand{\ev}{\text{\cyr \`{E}}}

\renewcommand{\a}{\alpha}

\def\G{\Gamma}
\def\b{\beta}

\def\vk{\varkappa}

\def\f{\varphi} 
\def\vk{\varkappa}
\def\s{\sigma}

 \DeclareMathOperator{\id}{id} 
  
\DeclareMathOperator{\Spec}{Spec}  
 \DeclareMathOperator{\im}{Im}

\DeclareMathOperator{\Hom}{Hom}  
 \DeclareMathOperator{\Sym}{Sym}

\DeclareMathOperator{\Diff}{Diff} 

\newcommand{\rw}{\rightarrow}
\newcommand{\lw}{\leftarrow}
\newcommand{\lrw}{\longrightarrow}
\newcommand{\llw}{\longleftarrow}

\newcommand{\Span}{\mathrm{span}}

\newcommand{\df}{\stackrel{\mathrm{def}}{=}}
\newcommand{\p}{\partial}

\newcommand{\ci}{C^{\infty}}
\newcommand{\Cc}{\mathcal C}
\newcommand{\Ee}{\mathcal E}
\newcommand{\Ff}{\mathcal F}
\newcommand{\Nn}{\mathcal N}
\newcommand{\Ll}{\mathcal L}

\newcommand{\Ss}{\mathcal S}

\newcommand{\Oo}{\mathcal O}
\newcommand{\Jj}{\mathcal J}

\newcommand{\Dd}{\mathcal D}
\newcommand{\Kk}{\mathcal K}
\newcommand{\Rr}{\mathcal R}

\newcommand{\dR}{\mathbb R}

\newcommand{\dC}{\mathbb C}

\newcommand{\gk}{\boldsymbol{k}}

\newcommand{\M}{\frak M}
\newcommand{\gS}{\frak S}

\DeclareFontFamily{OT1}{wncyi}{} \DeclareFontShape{OT1}{wncyi}{m}{it}{
   <5> <6> <7> <8> <9> gen * wncyi
   <10> <10.95> <12> <14.4> <17.28> <20.74> <24.88> wncyi10
  }{}
\DeclareSymbolFont{cyrletters}{OT1}{wncyi}{m}{it} 
\DeclareSymbolFontAlphabet{\cyrmath}{cyrletters} 
\DeclareMathSymbol{\rE}{\cyrmath}{cyrletters}{003} 
\DeclareMathSymbol{\rD}{\cyrmath}{cyrletters}{068} 
\DeclareMathSymbol{\rG}{\cyrmath}{cyrletters}{017} 
\DeclareMathSymbol{\rI}{\cyrmath}{cyrletters}{073} 
\DeclareMathSymbol{\rL}{\cyrmath}{cyrletters}{076} 
\DeclareMathSymbol{\rZ}{\cyrmath}{cyrletters}{090}

\newdimen\theight
\def \refright#1{%
             \vadjust{\setbox0=\hbox{\quad\vtop{\hsize5cm\bf\noindent #1}}%
             \theight=\ht0
             \advance\theight by \dp0    \advance\theight by \lineskip
             \kern -\theight \vbox to \theight{\rightline{\rlap{\box0}}%
             \vss}%
             }}%

\begin{document}
\title[What are symmetries]{\large   What are symmetries of nonlinear PDEs and what are they themselves? }\author{Alexandre Vinogradov} 
\date{\today}

\medskip

\abstract 
The general theory of (nonlinear) partial differential equations originated by S.~Lie had a significant
development in the past 30-40 years. Now this theory has solid foundations, a proper language,
proper techniques and problems, and a wide area of applications to physics, mechanics, to say nothing about traditional mathematics. However, the results of this development are not yet sufficiently
known to a wide public. An informal introduction in a historical perspective to this subject 
presented in
this paper aims to give to the reader an idea about this new area of mathematics and, possibly,
to attract new researchers to this, in our opinion, very promising area of modern mathematics.
\endabstract

\maketitle
\tableofcontents

This is neither a research nor a review paper but some reflections about general theory of 
(nonlinear) partial differential equations (N)PDEs and its strange marginal status in the realm 
of modern mathematical sciences.

For a long time a ``zoological-botanical" approach was and still continues to dominate in the study 
of PDEs, and, especially, of NPDEs. Namely, single equations coming from geometry, physics, mechanics, etc, were  ``tamed and cultivated" like single animals/plants of a practical or theoretical interest. As a rule, for each of these equations were found some prescriptions for the treatment motivated by some concrete {\it external}, i.e., physical, etc reasons, but  not based on the 
knowledge of their {\it intrinsic} mathematical nature.  Mainly, these prescriptions are focused
on how to construct  the solutions rather than to answer numerous questions  concerning global
properties of the PDE itself. 

Modern genetics explains what are living things, their variety and how to treat them to get the 
desired result. Obviously, a similar theory is indispensable for PDEs, i.e., a solid, well established 
{\it general theory}. The recent spectacular progress in genetics became possible only on the
basis of not less spectacular developments in chemistry and physics in the last century. Similarly, 
the general state of the art in mathematics 50-60 years ago was not  sufficiently mature to think 
about the general theory of PDEs. For instance, the fact that an advanced homological algebra 
will become an inherent feature of this theory could have been hard to imagine at that time.

Recent developments in general theory of PDEs are revealing more and more its intimate
relations with quantum mechanics, quantum field theory and related areas of contemporary
theoretical physics, which, also,  could be hardly expected a priori. Even more, now we can 
be certain that the difficulties and shortcomings of current physical theories are largely due to 
this historically explainable ignorance.

In this paper we informally present in a historical perspective problems, ideas and results that 
had led to  the renaissance of a general theory of PDEs after the long dead season that followed 
the pioneering S. Lie opera. Our guide was a modern interpretation of the Erlangen program
in the form of the principle : look for the symmetries and you will find the right way. Also, one 
of our goals 
was to show that this theory is not less noble part of pure mathematics than algebraic geometry,
which may be viewed as its zero-dimensional subcase. The paradox is that the number 
of mathematicians who worked on this theory does not exceed the number of those 
who studied Kummer surfaces. 

{\bf Warnings.} The modern general theory of PDEs is written in a new, not commonly known mathematical language, which was formed in the past 30-40 years and was used by a
very narrow circle of experts in this field. This makes it impossible to present this theory to 
a wide mathematical audience, to which this paper is addressed,  in its native language.
This is why the author was forced to be sometimes rather generic and to refer to some 
``common places" instead. His apology is in a maxim attributed to Confucius : ``An ordinary
man wonders marvelous things, a wise man wonders common places".

{\bf Notation.} Throughout the paper we use $\Lambda^k(M)$ (resp., $D_k(M)$) for the 
$\ci(M)$--module of $k$--th order differential forms (resp., $k$--vector fields) on a
smooth manifold $M$. For the rest the notation is standard.

\section{A brief history of nonlinear partial differential equations.}
Sophus Lie was the pioneer who sought for an order in the primordial chaos reigning in the 
world of NPDE at the end of the XIX century. The driving force of his approach was the idea to 
use symmetry considerations in the context of PDEs in the same manner as they were used 
by E.~Galois in the context of algebraic equations.  In the initial phase of realization of this 
program, Lie was guided by the principle ``chercher la sym\'etrie" and he discovered that behind numerous particular tricks found by hand in order to solve various concrete differential equations there are groups of transformations preserving  these equations, i.e., their symmetries.
Then, based on these ``experimental data", he developed the machinery of transformation
groups, which allows one to systematically compute what is now called {\it point} or {\it classical}
symmetries of differential equations.  Central in Lie theory is the concept of an {\it infinitesimal transformation} and hence of an {\it infinitesimal symmetry}. Infinitesimal symmetries of a 
differential equation or, more generally, of an object in differential geometry, form a Lie algebra.
This Lie's invention is among the most important in the history of mathematics.

Computation of classical symmetries of a system of differential equations leads to another 
{\it nonlinear} system, which is much more complicate than the original one. Lie resolved this
seemingly insuperable difficulty by passing to infinitesimal symmetries. In order to find them 
one has to solve an overdetermined system of {\it linear} differential equations, which is a much 
easier task and it non infrequently allows a complete solution. Moreover, by exponentiating 
infinitesimal symmetries one can find almost all finite symmetries. A particular case of this 
mechanism is the famous relation between Lie algebras and Lie groups.

Initial expectations that groups of classical symmetries are analogues of Galois group
for PDEs had led to a deep delusion. Indeed, computations show that this group for a 
generic PDE is trivial. This was one of the reasons why systematic applications of Lie 
theory to differential equations be frozen for a long time and the original intimate 
relations of this theory with differential equations were lost. Only much 
later in 1960-70  L.~V.~Ovsiannikov and his collaborators resumed these relations (see 
\cite{Ovs, Ibr}) and now they are extending in various directions.

Contact geometry was another important contribution of S. Lie to the general theory of NPDE's.
Namely, he discovered that symmetries of a first order NPDE imposed on one unknown 
function are {\it contact transformations}. These transformations not only mix dependent and 
independent variables but their derivatives as well. For this reason they are much more general 
than the above-mentioned point transformations, which mix only dependent and independent 
variables. Moreover, it turned out that contact symmetries are sufficient to build a complete
theory for this class of equations, which includes an elegant geometrical method of construction
of their solutions. In this sense contact symmetries play the role of a Galois group for this class 
of equations. On the other hand, the success of contact geometry in the theory of first order 
NPDEs led to the suspicion that classical symmetries form just a small part of all 
true symmetries of NPDEs. But the question of what are
these symmetries remained unanswered for a long time up to the discovery of {\it integrable 
systems} (see below). But some signs of an implicit use of such symmetries in some concrete situation can already be found in works of A.~V.~B\"{a}clund, E.~Noeter.

A courageous attempt to build a general theory of PDEs was undertaken by Charles Riquier 
at the very end of the XIX.th century.  ``Courageous" because at that time the only way to deal 
with general
PDEs was to manipulate their coordinate descriptions. His results then gathered in a handsome 
book  \cite{Riq} of more than 600 pages present, from the modern point of view, the first 
systematically developed general theory of formal integrability.  This book is full of cumbersome
computations, and the obtained results are mostly of a descriptive nature and do not reveal 
structural units of the theory. Nevertheless, it demonstrated that a general theory of 
PDE, even on a formal level, is not impossible. Moreover, Riquier showed that the formal theory 
duly combined with the Cauchy-Kowalewski theorem lead to various existence results in the
class of analytic functions such as the famous Cartan-K\"{a}hler theorem (see \cite{Cart, Kah}).  
In its turn Riquier's work motivated Elie Cartan  to look for a coordinate-free
language for  the formal/analytical theory and it led him to  the theory of {\it differential 
systems} based on the calculus of differential forms (see \cite{Cart}). Cartan's theory was 
later developed and extended by  E. K\"{a}hler \cite{Kah}, P.~K.~Rashevsky \cite{Rash}, 
Kuranishi \cite{Kura} and others. The reader will find its latest version in \cite{Bryant}.       

 In the middle of the XX.th century the  theory of  differential systems circulated  in a  narrow
 group of geometers as the most general theory of PDEs. However, this was an exaggeration.
For instance, there were no relations between this theory and the theory of linear PDEs, which 
was in a booming growth at that time. Moreover, this theory did not produce any, worth to 
be mentioned, application to the study of concrete PDEs. We can say that it is even hardly 
possible to imagine
that the study of the Einstein or Navier-Stokes equations will become easier after being converted
into differential systems.  So, the apparatus of differential forms  did not confirm the expectations
to become a natural base language for the general theory of PDEs, but it became one of the basic
instruments in modern differential geometry and in many areas of its applications.                                                     

The original  Riquier approach was improved and developed by Janet (\cite{Jan}). But,
unfortunately, his works were for a long period shadowed  by E.~Cartan works. Their vitality
was confirmed much later at the beginning of the new era for NPDEs (see \cite{Pom}). This 
era implicitly starts with the  concept of a jet bundle launched by Ch. Ehresmann (see \cite{Ehr}). Ehresmann himself did not develop applications of jet bundles to PDEs. But, fortunately,
this term became a matter of fashion and later was successfully used in various areas of differential
geometry.  In particular, D.~C.~Spencer and H.~Goldschmidt essentially used  jet bundles in 
their new theory of formal integrability by inventing a new powerful instrument, namely, the 
Spencer cohomology (see \cite{Spen, Pom}). In this way were discovered the first structural 
blocks of the general theory of PDEs, and this new theory demonstrated some important 
advantages in comparison with the theory of differential systems.                   

The Ehresmann concept of jet bundle is, however, too restrictive to be applied to general
NPDEs and for this reason should be extended to that of jet  space (or manifold). Namely, 
the $k$--th 
order jet bundle $J^k(\pi)$ associated with a smooth bundle $\pi:E\rightarrow M$ consists of
$k$--th order jets of sections of $\pi$, while the $k$--th order jet space $J^k(E,n)$ consists of
$k$--th order jets of $n$--dimensional submanifolds of the manifold $E$. Jet spaces are naturally supplied with a structure, called the {\it Cartan distribution}, which allows an interpretation of 
functions on 
them as nonlinear differential operators. Differential equations in the standard but coordinate 
independent meaning of this term are naturally interpreted as submanifolds of jet spaces. 
The first systematic study of geometry of jet spaces was done by A.~M.~Vinogradov and 
participants of his Moscow seminar in the seventieth  of the passed century 
(see \cite{V1,Klava,KLV}). Later, on this basis, it was understood that various
natural differential operators and constructions that are necessary for the study of a system of
PDEs of  order $k$ do not live necessarily on the $k$--th order jet space but involve jet spaces 
of any order.
This is equivalent to say that a conceptually complete theory of PDEs  is possible only on {\it  
infinite order} jet spaces.  A logical consequence of this fact is that objects of the {\it category of
partial differential equations} are {\it diffieties}, which duly formalize the vague idea of the ``space 
of all solutions" of a PDE. Diffieties are a kind of infinite dimensional manifolds, and the specific differential calculus on them, called {\it secondary calculus}, is a native language to deal with
PDEs and especially with NPDEs (see \cite{V2,KVe,KVin}).

Below we shall show how to come to secondary calculus and hence to the {\it general theory}
of (nonlinear) partial differential equations by trying to answer the question ``what are symmetries
of a PDE". It is worth stressing that Klein's Erlangen program was a good guide in this expedition, which was decisive in finding the right way in some crucial moments.

\section{Evolution of the notion of symmetry for differential equations.}
A retrospective view on how the answer to the question ``what are symmetries of a PDE" 
evolved  historically will be instructive for our further discussion. From the very beginning 
this question was more implicitly than explicitly related with the answer to question ``what 
is a PDE". It seems that the apparent absurdity of this question prevented its exact 
formulation and hence slowed the development of the general theory.

Below we shall use the following notation. If $\sigma=(\sigma_1,\dots,\sigma_n)$ is a multiindex,
then $|\sigma |=\sigma_1+\dots+\sigma_n$, and
$$
\frac{\partial^{|\sigma |} f(x)}{\partial x^\sigma}=
\frac{\partial^s f(x)}{\partial x_1^{\sigma_1}\dots x_n^{\sigma_n}},\quad |\sigma |=s, \quad
f(x)=f(x_1,\dots,x_n).
$$
We assume that  $\frac{\partial^{|\sigma |} f}{\partial x^\sigma}=f$ if $\sigma=(0,\dots,0)$.
\newline

{\bf Standard (``classical") definition.} According to the commonly accepted point of view a 
system of PDEs is a set of expressions
\begin{equation}\label{eq1.1}
F_i(x, u, u_{[1]},\dots, u_{[k]})=0,\quad i=1, 2,\dots, l.
\end{equation}
where $x=(x_1,\dots, x_n)$ are {\it independent variables}, $u=(u^1,\dots, u^m)$  {\it dependent} 
ones, and $u_{[s]}$ stands for the totality of symbols $u^i_{\sigma}, \,1\le i \le m,$ with $|\sigma |=s$.
Further we shall use short ``PDE" for ``system of PDEs" and, accordingly, write 
$$
F(x, u, u_{[1]},\dots, u_{[k]})=0 \quad \mbox{assuming that}\quad F=(F_1,\dots, F_l).
$$
{\bf Solutions.} A system of functions $f_1(x),\dots,f_m(x)$ is a solution of the PDE (\ref{eq1.1}) 
if the substitutions $\frac{\partial^{|\sigma |} f^i}{\partial x^\sigma}\to u^i_{\sigma}$ transform the expressions (\ref{eq1.1}) to identically equal to zero functions of $x$.

This traditional view on PDEs is presented in all, modern and classical, textbooks. For instance, in 
Wikipedia one may read that a PDE is ``an equation that contains unknown multivariable functions 
and their partial derivatives" or ``une \'equation aux d\'erivŽes partielles (EDP) est une \'equation 
dont les solutions sont les fonctions inconnues v\'erifiant certaines conditions concernant leurs 
d\'eriv\'es partielles."\newline

{\bf Symmetries: the first idea.} The  ``common sense"  coherent with this point of view
suggests to call a symmetry of PDE (\ref{eq1.1}) transformations 
\begin{equation}\label{tr1}
x_i=\phi_i(\bar{x}_1,\dots,\bar{x}_n), \,i=1,\dots,n, \quad 
u^j=\psi^j(\bar{u}^1,\dots,\bar{u}^m), \,j=1,\dots,m,
\end{equation} 
of dependent and independent variables that ``preserve the form" of relations (\ref{eq1.1}). 
More exactly, this means that the so-obtained functions 
$\bar{F_i}=\bar{F_i}(\bar{x},\bar{u},\bar{u}_{[1]},\dots,\bar{u}_{[s]})$s are linear combinations 
of functions $F_i$s with functions of $x, u, u_{[1]},\dots, u_{[k]}$ as coefficients. Here we  
used the confusing classical notation where $(x,u)$ stands for coordinates of the image of the 
point $(\bar{x},\bar{u})$. Also, it is assumed that transformations of symbols $u_{\sigma}^i$ 
are naturally induced by those of $x$ and $u$.

Many fundamental equations in physics and mechanics inherits space-time 
symmetries, and these are ``first idea" symmetries. Very popular in mechanics of continua 
dimensional analysis is also based on the so-understood concept of 
symmetry (see \cite{Bar, BluC, Ovs}).
\begin{ex}
The Burgers equation $u_t=u_{xx}+uu_x$ is invariant, i.e., symmetric, with respect to space
shifts $(x=\bar{x}+c, t=\bar{t}, u=\bar{u})$, time shifts $(x=\bar{x}, t=\bar{t+c}, u=\bar{u})$ and
the passage to another Galilean inertial frame moving with the velocity $v$ 
$(x=\bar{x}+vt, t=\bar{t}, u=\bar{u})$. This equation possesses also {\it scale} symmetries:
$x=\lambda\bar{x}, \,t=\lambda^2\bar{t}, \,u=\lambda^{-1}\bar{u}, \,\lambda\in\dR$.
\end{ex}

The above definition of a symmetry is based on the a priori premise that the division of variables 
into dependent and independent ones is an indispensable  part of the definition of a PDE. 
However, many arguments show that this point of view is too restrictive. In particular, numerous
tricks that were found by hands to resolve various concrete PDEs involves transformations which 
do not respect this division. For instance, transformations 
$$
x=\frac{\bar{x}}{\tau\bar{t}+1}, \quad t=\frac{\bar{t}}{\tau\bar{t}+1},  \quad
u=\bar{u}+\tau(\bar{t}\bar{u}-\bar{x})
$$
depending on a parameter $\tau\in\dR$ leave the Burgers equations invariant. They, however, 
do not respect sovereignty of of the dependent variable $u$. 

In this connection more obvious and important argument is that 
\begin{quote}
{\it what is called functions in the traditional definition of a PDE
are not, generally, functions but elements of coordinate-wise descriptions of certain objects, like tensors, submanifolds, etc.}
\end{quote} 
Indeed, if dependent variables $u$ are components of a tensor, then a transformation of
independent variables induces automatically a transformation of independent ones. So, the 
division of variables into dependent and independent ones cannot, in principle, be respected 
in such cases. Moreover, this, as banal as well-known observation, which is nevertheless 
commonly ignored, poses a question
\begin{quote}
{\it what are ``independent variables", i.e., what are mathematical objects that are subjected by 
PDEs?}
\end{quote} 
The ``obvious" answer that these are ``objects that are described coordinate-wisely by means of functions" is purely descriptive and hence not very satisfactory. In fact, this question is neither trivial, nor stupid, and, in particular, its analysis directly leads to the conception of jets (see below). \newline

{\bf Symmetries: the second idea.}
Under the pressure of the the above arguments it seems natural to call a symmetry of PDE (\ref{eq1.1})
a transformation of independent and dependent, which respect the status of independent 
variables  only, i.e., a transformation of the form
\begin{equation}\label{tr2}
\left\{
\begin{array}{lll}
x_i=&\phi_i(\bar{x}_1,\dots,\bar{x}_n), \quad i=1,\dots,n, \quad \\
u^j=&\psi^j(\bar{x}_1,\dots,\bar{x}_n,\bar{u}^1,\dots,\bar{u}^m), \quad j=1,\dots,m.\\
\end{array}
\right.
\end{equation}
This idea is consistent with many situations in physics and mechanics where space-time coordinates play role of independent variables, while `internal" characteristics of the considered
continua, fields, etc, refer to dependent ones. Mathematically, these quantities are represented as
sections of suitable fiber bundles, and transformations that preserve the bundle structure are
exactly of the form (\ref{tr2}). Gauge transformations in modern physics are of this kind.

On the other hand,  since the second half of 18th century, the development of differential 
geometry put in light various problems related with 
surfaces and, later, with manifolds and their maps (see \cite{Mng,Gau,Darb}) formulated in  
terms of PDEs. A surface in the 3-dimensional Euclidean space $E^3$ is not, generally, the 
graph of a function. So, the phrase that the equation 
\begin{equation}\label{min}
(1+u_x^2)u_{yy}-2u_xu_yu_{xy}+(1+u_y^2)u_{xx}=0
\end{equation}
is the equation of minimal  surfaces is not, rigorously speaking, true. More exactly, it is true
only locally for surfaces of the form $z=u(x,y)$ with $(x,y,z)$ being  standard Cartesian 
coordinates in $E^3$. So, the question ``what is the true (global) equation of minimal  
surfaces" should be clarified. This
question, which was historically ignored, becomes, however, rather relevant if one thinks about
global topological properties of minimal surfaces. Also, isometries of $E^3$ preserve the class
of minimal surfaces and hence they must be considered as symmetries of the ``true" equation of minimal  surfaces in any reasonable sense 
of this term. But, generally, these transformations do not respect 
the status of both independent and dependent variables. This and many other similar
examples show that the second idea is still too restrictive.

{\bf Symmetries: the third idea.} The next obvious step is to consider transformations
\begin{equation}\label{tr2}
\begin{array}{rcr}
x_i=\phi_i(\bar{x}_1,\dots,\bar{x}_n,\bar{u}^1,\dots,\bar{u}^m), \quad
&u^j=\psi^j(\bar{x}_1,\dots,\bar{x}_n,\bar{u}^1,\dots,\bar{u}^m),\\
&i=1,\dots,n,  \quad j=1,\dots,m.\\
\end{array}
\end{equation}
as eventual symmetries of PDEs. In this form the idea of a symmetry of a PDE was formulated 
by S.~Lie at the end of 19th century and was commonly accepted for a long time up to the 
discovery of integrable systems at the late 1960s, when some strong doubts about it arose. 
Symmetries of form (\ref{tr2}) are called {\it point symmetries} in order to distinguish them
from {\it contact symmetries} (see below) and more general ones that recently emerged.

{\bf Symmetries: the fourth incomplete idea.} But S.~Lie
himself created the ground for such doubts by developing the theory of one first order PDEs in
the form of contact geometry. From this point of view natural candidates for symmetries of such
a PDE are contact transformations, which mix independent and dependent variables and their 
first order derivatives in an almost arbitrary manner. In particular, this means that, generally, transformations of dependent and independent variables involve also first derivatives. i.e.,
\begin{equation}\label{tr3}
\left\{
\begin{array}{lll}
x_i=&\phi_i(\bar{x}_1,\dots,\bar{x}_n,\bar{u},\bar{u}_{x_1},\dots,\bar{u}_{x_m}), \quad i=1,\dots,n, \quad \\
u=&\psi^j(\bar{x}_1,\dots,\bar{x}_n,\bar{u},\bar{u}_{x_1},\dots,\bar{u}_{x_m}).\\
\end{array}
\right.
\end{equation}
Transformations (\ref{tr3}) are to be completed by transformations of first derivatives
$$
x_i=\phi_i(\bar{x}_1,\dots,\bar{x}_n,\bar{u},\bar{u}_{x_1},\dots,\bar{u}_{x_m}), \quad i=1,\dots,n,
$$
in a way that respects the ``contact condition" $d\bar{u}-\sum_{i=1}^n\bar{u}_{x_i}dx_i=0$.

Contact transformations can be naturally prolonged to transformations of higher order derivatives
and, therefore, considered as candidates for true symmetries of PDEs with one dependent variable.  For instance, as such they are very useful in the study of Monge-Amp\'ere equations (see \cite{KLR}). In other words, the third idea becomes too restrictive,  at least, for equations
with one dependent variable.

The above discussion leads to a series of questions: 

QUESTION 1. What are analogues of contact transformations for PDEs with more than one dependent variable? 

QUESTION 2. Are there higher order analogues of contact transformations, i.e., transformations
mixing dependent and independent variables with derivatives of order higher than one?

To answer this questions we, first, need to bring the traditional approach to PDEs to a more
conceptual form. In particular, a coordinate-free definition of a PDE equivalent to the standard 
one is needed. This is done in the next section.

\section{Jets and PDEs.}
Various objects (functions, tensors, submanifolds, smooth maps, geometrical structures, etc.) 
that are subjected to PDEs may be interpreted as submanifolds of a suitable manifold. For
instance, functions, sections of fiber bundles, in particular, tensors, and smooth maps may be geometrically viewed as the corresponding graphs. So, we  assume this unifying point of 
view and interpret PDEs as {\it differential restrictions} imposed on submanifolds of a given
manifolds.

\subsection{Jet spaces.}\label{JS}
So, objects of our further considerations will be $n$--dimensional submanifolds of an 
$(n+m)$--dimensional submanifold $E$.
Let $L\subset E$ be a such one. In order to locally describe it  in a local chart 
$(y_1,\dots,y_{n+})$ we must choose among these coordinates $n$ independent on $L$ ones, 
say, $(y_{i_1},\dots,y_{i_n})$ and declare the remaining $y_j$s to be dependent ones. The 
notation $x_1=y_{i_1},\dots,x_n=y_{i_n}, \,u^1=y_{j_1}\dots,u^m=y_{j_m}$ with 
$\{j_1,\dots,j_m\}=\{1,\dots,n+m\}\setminus\{i_1,\dots,i_n\}$ stresses this artificial division 
of local coordinates into dependent and independent ones. We shall refer to $(x,u)$ as a 
{\it divided chart}. By construction $L$ is locally described in this divided chart by equations 
of the form $u^i=f^i(x)$, $i=1,\dots,n$. The next step is to understand what is the manifold in 
which $(x,u,u_{[1]}\dots,u_{[k]})$ is a local chart. The answer is as follows.

Let $M$ be a manifold, $z\in M$ and $\mu_z=\{f\in C^{\infty}(M)\,|\,f(z)=0\}$ the ideal of the 
point $z$. Elements of the quotient algebra $J_z(M)=C^{\infty}(M)/\mu_z^{k+1}$ ($\mu_z^s$ 
stands for the $s$-th power of the ideal $\mu_z$) are called {\it $k$--th order jets} of functions 
at the point $z\in M$. The {\it $k$--th jet of $f$ at $z$} denoted by $[f]_z^k$ is the image of $f$ 
under the factorization homomorphism $C^{\infty}(M)\to J_z(M)$. This definition also holds for 
$k=\infty$ if we put $\mu_z^{\infty}=\bigcap_{k\in \mathbb{N}}\mu_z^k$. It is easy to see that 
$[f]_z^k$=$[g]_z^k$ if and only if in a local chart all the derivatives of the functions $f$ and 
$g$ of order $\leq k$ at the point $z$ are equal.

Two $n$--dimensional submanifolds $L_1,L_2\subset E$ are called {\it tangent with the order 
$k$ at a common point $z$} if for any $f\in C^{\infty}(M) \;[f|_{L_1}]_z^k=0$ implies 
$[f|_{L_2}]_z^k=0$ and vice versa. Obviously, $k$-th order tangency is an equivalence relation.
\begin{defi}\label{jet}
The equivalence class of $n$--dimensional submanifolds of $E$, which are $k$-th order tangent  
to $L$ at $z\in L$,  is called the \emph{$k$-th order jet of $L$ at $z$} and is denoted by $[L]_z^k$.
\end{defi}
The set of all $k$-jets of $n$--dimensional submanifolds $L$ of $E$ is naturally supplied
with the structure of a smooth manifold, which will be denoted by $J^k(E,n)$. Namely, associate
with an $n$--dimensional submanifold $L$ of $E$ the map
$$
j_k(L):L\rightarrow J^k(E,n), \,L\ni z\mapsto [L]_z^k.
$$
and call a function $\phi$ on $J^k(E,n)$ {\it smooth} if $j_k(L)^*(\phi)\in C^{\infty}(L)$ for all
$n$--dimensional $L\subset E$. The so-defined smooth function algebra will be denoted by
$\Ff _k(E,n)$, i.e., $C^{\infty}(J^k(E,n))=\Ff_k (E,n)$. 
\begin{rmk}
If $k<\infty$ the above definition of the smooth structure on $J^k(E,n)$ is equivalent to the 
standard one, which use charts and atlases (see below). But it becomes essential for $k=\infty$, 
since the standard ``cartographical" approach in this case creates some boring inconveniences.
\end{rmk}
If $L$ is given by equations $u^i=f^i(x)$, $i=1,\dots,n$, in a divided chart and
$(x_1^0,\dots,x_n^0,u_0^1\dots,u_0^m)$ are coordinates of $z$ in this chart, then, as it is 
easy to see, $[L]_z^k$ is uniquely defined by the derivatives
\begin{equation}\label{j-chart}
u_{\sigma,0}^i=\frac{\partial^{|\sigma |} f^i}{\partial x_\sigma}(x_1^0,\dots,x_n^0),\quad
1\le i\le m,\quad |\sigma|\le k,
\end{equation}
and vice versa. So, the numbers $x_j^0$ together with the numbers $u_{\sigma,0}^i$  may 
be taken for local coordinates of the point $\theta=[L]_z^k\in J^k(E,n)$. By observing that
$$
u_{\sigma,0}^i=j_k(L)^*\left(\frac{\partial^{|\sigma |} f^i}{\partial x_\sigma}\right)(\theta)
$$
we see that the functions  
$u_{\sigma}^i\df j_k(L)^*\left(\frac{\partial^{|\sigma |} f^i}{\partial x_\sigma}\right), 
\,1\leq i\leq m, \,|\sigma|\leq k,$m together with the functions $x_j$ form a smooth local 
chart on $J^k(E,n)$. 

Thus we see that $(x,u,u_{[1]},\dots,u_{[k]})$ is a local chart on $J^k(E,n)$ and hence (\ref{eq1.1}) 
is the equation of a submanifold in $J^k(E,n)$. This allows us to interpret the standard definition
of PDEs in an invariant coordinate-free manner.
\begin{defi}\label{eq.jet}
A system of PDEs of order $k$ imposed on $n$--dimensional submanifolds of a manifold 
$E$ is a submanifold $\Ee$ of $J^k(E,n)$.
\end{defi}
\begin{rmk}
$\Ee$ as a submanifold of $J^k(E,n)$ may have singularities.
\end{rmk}

\subsection{Jet tower.} 
Note that $E$ is naturally identified with 
$J^0(E,n):\,z\leftrightarrow [L]_z^0$, and natural projections
$$
\pi_{k,l}:J^k(E,n)\rightarrow J^l(E,n), \;[L]_z^k\mapsto [L]_z^l, \,l\leq k,
$$ 
relate jet spaces of various orders in a unique structure
\begin{equation}\label{j-tower}
E=J^0(E,n)\stackrel{\pi_{1,0}}{\longleftarrow} J^1(E,n)\stackrel{\pi_{2,1}}{\longleftarrow}\dots
\stackrel{\pi_{k,k-1}}{\longleftarrow}J^k(E,n)\stackrel{\pi_{k+1,k}}{\longleftarrow}\dots J^{\infty}(E,n).
\end{equation}
It is easy to see that  $J^{\infty}(E,n)$ is the inverse limit of the system of maps $\{\pi_{k,l}\}$.
Also note that $\pi_{k,l}:J^k(E,n)\to J^l(E,n)$ is a fiber bundle. Moreover, 
$\pi_{k,k-1}:J^k(E,n)\to J^{k-1}(E,n)$ is an {\it affine bundle} if $k\geq 2$ and $m>1$ or if
$k\geq 3$ and $m=1$ (see \cite{KLV,Pom}).

Dually to (\ref{j-tower}), smooth function algebras on jet spaces form a telescopic system of
inclusions
\begin{equation}\label{F-tower}
\ci(E)=\Ff\stackrel{\pi_{1,0}^*}{\rightarrow} \Ff_1\stackrel{\pi_{2,1}^*}{\longrightarrow}\dots
\stackrel{\pi_{k,k-1}^*}{\longrightarrow}\Ff_k\stackrel{\pi_{k+1,k}^*}{\longrightarrow}
\dots \Ff_{\infty}.
\end{equation}
So, $\Ff_{\infty}$ may be viewed as the direct limit of (\ref{F-tower}). By identifying $\Ff_k$ with 
$\pi_{\infty,k}^*(\Ff_k)$ we get the filtered algebra 
$\Ff_0\subset\Ff_1\subset\dots\subset\Ff_k\subset\dots, \;\Ff_{\infty}=
\bigcup_{k=0}^{\infty}\Ff_k$.

Since two submanifolds of the same dimension are first order tangent at a point $z$ if and
only if they have the common tangent space at $z$, $[L]_z^1$ is naturally identified with 
$T_zL$. In this way
$J^1(E,n)$ is identified with the Grassmann bundle $\textrm{Gr}_n(E)$ of n-dimensional 
subspaces tangent to E , and the canonical projection $\textrm{Gr}_n(E)\to E$ is identified 
with $\pi_{1,0}$. In particular, the standard fiber of $\pi_{1,0}$ is the Grassmann manifold 
$Gr_{n,m}$so that $\pi_{1,0}$ is not an affine bundle. If $m=1$, then the fiber of $\pi_{2,1}$ 
is the Lagrangian Grassmanian.
\begin{ex}
The equation $\Ee$ of minimal surfaces in the 3-dimensional Euclidean space $E^3$ is a 
hypersurface in $J^2(E^3,2)$. The projection $\pi_{2,1}:\Ee\to J^1(E^3,2)$ is a nontrivial
bundle whose fiber is the 2-dimensional torus. So, rigorously speaking, \emph{(\ref{min})} 
is not the equation of minimal surfaces but a local piece of it. 
\end{ex}
 This and many other similar examples show that, generally, (\ref{eq1.1}) is just a local 
 coordinate-wise description of a PDE.
 
 \subsection{Classical symmetries of PDEs.}\label{cl-sy}
 The language of jet spaces not only gives a due conceptual rigor to the traditional theory of
 PDEs but it also simplifies many technical aspects of it and makes transparent and better workable
 various basic constructions. This will be shown in the course of the subsequent exposition.
 But now we shall illustrate this point by explaining how ``point transformations" acts on PDEs.
 
 First, we observe that (\ref{tr2}) is just a local coordinate-wise description of a diffeomorphism
 $F:E\to E$. Now the question we are interested in is: how does $F$ act on jets. The answer is obvious: $F$ induces the diffeomorphism 
 \begin{equation}
F_{(k)}:J^k(E,n)\to J^k(E,n), \;[L]_z^k\mapsto [F(L)]_{F(z)}^k,
\end{equation}
called the {$k$--lift} of $F$. This immediately leads to formulate the definition of a ``classical" (=``point") symmetry of a PDE.
\begin{defi}\label{cl-sym}
A classical/point symmetry of a PDE $\Ee\subset J^k(E,n)$ is a diffeomorphism $F:E\to E$ 
such that $F_{(k)}(\Ee)= \Ee$.
\end{defi}
Similarly one can define lifts of ``infinitesimal point transformations", i.e., vector fields on $E$. 
Recall that if $X$ is a vector field on $E$ and $F_t:E\to E$ is the flow generated by it , then 
$$
X=\left.\frac{d(F_t^*)}{dt}\right |_{t=0} \mbox{with} \;F_t^*:C^{\infty}(E)\to C^{\infty}(E).
$$
Then the lift $X_{(k)}$ of $X$ to $J^l(E,n)$ is defined as
$$
X_{(k)}=\left.\frac{d((F_t)_{(k)}^*)}{dt}\right |_{t=0} \;.
$$
\begin{defi}
A vector field $X$ on $E$ is an \emph{infinitesimal classical/point symmetry of a PDE $\Ee\in J^k(E,n)$} if $X_{(k)}$ is tangent to $\Ee$.
\end{defi}
Nonlinear partial differential operators are also easily defined in terms of jets.
\begin{defi}
A \emph{(nonlinear) partial differential operator}  $\Box$ of order $k$ sending $n$--dimensional 
submanifolds of $E$ to $E'$ is defined as the composition $\Phi\circ j_k$ with 
$\Phi:J^k(E,n)\to E'$ being a (smooth) map, i.e.,  $\Box(L)=\Phi(j_k(L))$.
\end{defi}
\noindent For instance, functions on $J^l(E,n)$ are naturally interpreted as (nonlinear)
differential operators.
 
The above definitions and constructions are easily specified to fiber bundles. Namely, if
$\pi:E\to M, \,\dim\,M=n$, is a fiber bundle, then the $k$--th order jet $[s]_x^k$ of an local 
section of it $s:U\to E$ ($U$ is an open in $M$) at a point $x\in M$ is defined as 
$[s]_x^k=[s(U)]_{s(x)}^k$. These specific jets form an everywhere dense open subset in
$J^k(E,n)$ denoted by $J^k(\pi)$ and called the {\it $k$--order jet bundle of $\pi$}. The 
substitute of maps $j(L)$ in this
context are maps $j_k(s)\df j_k(s(U))$. Additionally, we have naturaI projections 
$\pi_k=\pi\circ \pi_{k,0}:J^k(\pi)\to M$. If $\pi$ is a vector bundle, then $\pi_k$ is a vector bundle
too, and the equation $\Ee\subset J^k(\pi)$ is {\it linear} if $\Ee$ is a linear sub-bundle of $\pi_k$,
etc.  For further details concerning the ``fibered" case, see  \cite{KVin, Klava, Pom}.

\section{Higher order contact structures and generalized solutions of NPDEs.}
\label{high-cont-str}
\subsection{Higher order contact structures.}
Now we are going to reformulate the standard definition of a solution of a PDE in a 
coordinate-free manner. Put $L_{(k)}=\im j_k(L)$ for an $n$--dimensional submanifold 
of $E$. Obviously, $L_{(k)}$ is an $n$--dimensional submanifold $L$ of $J^k(E,n)$, which 
is projected diffeomorphically onto $L$ via $\pi_{k,0}$.
\begin{defi}\label{j-sol}
L is a \emph{solution  in the standard sense of a PDE $\Ee\subset J^k(E,n)$} 
if $L_{(k)}\subset \Ee$.
\end{defi}

If $u^i=F^i(x), i=1,\dots,n$ are local equations of $L$, then
$$
u_{\sigma}^i=\frac{\p^{|\sigma|}f^i}{\p x^{\sigma}}(x), \quad i=1,\dots,n, \quad |\sigma|\leq n,
$$
are local equations of $L_{(k)}$ in $J^k(E,n)$. This shows that the coordinate-free definition
\ref{j-sol} coincides with the standard one. Also, we see that the $L_{(k)}$'s form a very special
class of $n$--dimensional submanifolds in $J^k(E,n)$. This class is not intrinsically defined,
and hence Definition \ref{j-sol} is not {\it intrinsic}. For this reason it is necessary to supply
$J^k(E,n)$ with an additional structure, which allows to distinguish the submanifolds 
$L_{(k)}$ from others. Such a structure is a distribution on $J^k(E,n)$ defined as follows.
\begin{defi}
The minimal distribution $\Cc^k:J^k(E,n)\ni\theta\mapsto\Cc_\theta^k\subset T_\theta(J^k(E,n))$
on $J^k(E,n)$ such that  all the $L_{(k)}$ are integral submanifolds of it, i.e., 
$T_{\theta}(L_{(k)})\subset\Cc_\theta^k, \,\forall\theta\in L_{(k)}$,  is called the 
\emph{ $k$--th order contact structure} or the \emph{Cartan distribution} on $J^k(E,n)$. 
\end{defi}
It directly follows from the definition that
\begin{equation}\label{span}
\Cc_\theta^k= \Span\{T_{\theta}(L_{(k)}) \;\mbox{for all} \;L \;\mbox{such that} \;L_{(k)}\ni\theta\}.
\end{equation}
\noindent Due to the importance of the subspaces $T_{\theta}(L_{(k)})\subset T_{\theta}(J^k(E,n))$ 
we shall call them {\it $R$--planes} (at $\theta$). By construction any $R$--plane at $\theta$ belongs
to $\Cc_{\theta}^k$. The following simple fact is very important and will 
be used in various constructions further on.
\begin{lem}\label{R_theta}
Let $\theta=[L]_z^k=[N]_z^k$ and $\theta '=\pi_{k,k-1}(\theta)$. Then 
$T_{\theta'}(L_{(k-1)})=T_{\theta'}(N_{(k-1)})$ and hence the $R$--plane
$R_{\theta}=T_{\theta'}(L_{(k-1)})$ is uniquely defined by $\theta$. Moreover,
the correspondence $\theta\mapsto R_{\theta}$ between points of $J^{k}(E,n)$ and 
$R$--planes at points of $J^{k-1}(E,n)$ is biunique.
\end{lem}
\noindent This lemma allows to identify the fiber $\pi_{k,k-1}^{-1}(\theta')$ with the variety of 
all $R$--planes at $\theta$ and hence  $J^k(E,n)$ with the variety of $R$--planes at
points of $J^{k-1}(E,n)$.

Below we list some basic facts concerning of the Cartan distribution and R-planes 
(see \cite{V1, KLV,Klava}).
\begin{proc}\label{C-list}
\begin{enumerate}
\item $\Cc_\theta^k=(d_{\theta}\pi_{k,k-1})^{-1}(R_{\theta})$ with 
$d_{\theta}\pi_{k,k-1}:T_{\theta}(J^{k}(E,n))\to T_{\theta'}(J^{k-1}(E,n))$ being the differential of
$\pi_{k,k-1}$ at $\theta$. In particular, \\$d_{\theta}\pi_{k,k-1}(\Cc_{\theta}^k) \subset \Cc_{\theta'}^{k-1}$.
\item In local coordinates the Cartan distribution is given by the equations
$$
\omega_{\sigma}^i\df du_\sigma^i-\sum_{j} u^i_{\sigma+1_j}dx_j=0, \quad |\sigma |<k, 
\;\mbox{where} 
\;(\sigma+1_j)_i=\sigma_i+\delta_{ij},
$$
or, dually, is generated by the vector fields 
$$
D_i^k=\frac{\p}{\p x_i}+\sum_{j, |\sigma|<k}  u_{\sigma+1_i}^j\frac{\p}{\p u_{\sigma}^j},
\quad i=1,\dots,n, \;\mbox{and} \;\frac{\p}{\p u_{\sigma}^j}, \;|\sigma|=k.
$$
\item
$$
\dim\Cc_\theta^k=m\left(\begin{array}{c}n+k-1 \\k\end{array}\right)+n, \;\mbox{if} \quad 0\leq k<\infty;
\quad \dim\Cc_\theta^{\infty}=n.
$$
\item Tautologically, a point $\theta=[L]_z^{\infty}\in J^{\infty}(E,n)$ is the inverse limit of
$\theta_k=\pi_{\infty,k}(\theta)=[L]_z^k, \,k\to\infty$. Then $\Cc_{\theta}^{\infty}$ is the inverse 
limit of the chain 
$$
\dots\stackrel{d_{\theta_k}\pi_{k,k-1}}{\longleftarrow}\Cc_{\theta_k}^k
\stackrel{d_{\theta_{k+1}}\pi_{k+1,k}}{\longleftarrow}\Cc_{\theta_{k+1}}^{k+1}
\stackrel{d_{\theta_{k+2}}\pi_{k+2,k+1}}{\longleftarrow}\dots
$$
\item Distributions $\Cc^k, \,k<\infty$, are, in a sense, ``completely non-integrable", while
their inverse limit $\Cc^{\infty}$ is completely (Frobenius) integrable and locally generated
by commuting {\it total derivatives}  
$$
D_i=\frac{\p}{\p x_i}+\sum_{j, \sigma}  u_{\sigma+1_i}^j\frac{\p}{\p u_{\sigma}^j},
\quad i=1,\dots,n.
$$ 
\item If an $n$--dimensional integral submanifold $N$ of $\Cc^k, \,k<\infty$, is transversal to 
fibers of $\pi_{k,k-1}$, then, locally, $N$ is of the form $L_{(k)}$ and, therefore, $\pi_{k,0}(N)$
is an {\it immersed} $n$--dimensional submanifold of $E$.
\end{enumerate}
\end{proc}

{\it Cartan's forms} $\omega_{\sigma}^i$ figuring in assertion (2) of the above proposition were systematically used by E.~Cartan in his reduction of PDEs to exterior differential systems. Hence 
the term ``Cartan distribution".

Note that if $m=1$, then the manifold $J^1(E,n)$ supplied with the Cartan distribution $\Cc^1$ 
is a contact manifold. The contact distribution $\Cc^1$ is locally given by the classical contact 
form $du-\sum_{i=1}^{n}u^idx_i=0$. So, $\Cc^k$ whose construction word for word mimics 
the classical construction of contact geometry may be viewed as its higher order analogue,
i.e.,  the $k$--th order contact structures. 

Recall now how the theory of one 1-st order PDE with one independent variable is formulated 
in terms of contact geometry. Let $K, \dim\,K=r+1$, be a manifold supplied with an  
$r$--dimensional distribution $\Cc:K\ni x\mapsto \Cc_x\subset T_xK$. The fiber at $x\in K$
of the {\it normal} to the $\Cc$ vector bundle $\nu_{\Cc}:N_{\Cc}\rw K$  is $T_xK/\Cc_x$, and
$\dim\,\nu_{\Cc}=1$. We shall write $X\in \Cc$ if the vector field $X$ belongs to $\Cc$, i.e.,
$X_x\in\Cc_x, \,\forall x\in K$. By abusing language we shall denote also by $\Cc$ the
$\ci(K)$--module of vector fields belonging to $\Cc$  and put $\Nn_{\Cc}=\Gamma(\nu_{\Cc})$.
The {\it curvature} of $\Cc$ is the following $\ci(K)$--bilinear skew-symmetric form $\Omega$ with
values in $\Nn_{\Cc}$:
$$
\Omega_{\Cc}(X,Y)=[X,Y] \mod \Cc, \quad X,Y\in \Cc.
$$
\noindent $\Omega_{\Cc}$ is {\it nondegenerate} if the map 
$$
\Cc\ni X\mapsto \Omega_{\Cc}(X,\cdot)\in\Lambda^1(K)\otimes_{\ci(K)}\Nn_{\Cc}
$$
is an isomorphism of $\ci(K)$--modules. The pair $(K,\Cc)$ is a {\it contact manifold} if the
2-form $\Omega_{\Cc}$ is nondegenerate. In such a case $r$ is odd, say, $r=2n+1$.

This definition of contact manifolds is not standard (see \cite{Arn,KLR}) but is more convenient 
for our purposes. By the classical Darboux lemma a contact manifold locally possesses 
{\it canonical coordinates} $(x_1,\dots,x_n,u,p_1,\dots,p_n)$ in which $\Cc$ is given by
the 1-form $\omega\df du-\sum_{i=1}^{n}p_idx_i=0$. Then $\boldsymbol{e}=[\p/\p u \mod \Cc]$
is a local base of $\Nn_{\Cc}$ and 
$\Omega_{\Cc}=-d\omega\otimes\boldsymbol{e}=
(\sum_{i=1}^{n}dp_i\wedge dx_i)\otimes\boldsymbol{e}$.

If a hypersurface $\Ee\subset K$ is interpreted as a 1-st order PDE, then a (generalized)
solution of $\Ee$ is a $Legendrian$ submanifold $L$  in $K$ belonging to $\Ee$. Recall that a Legendrian submanifold $L$ is an $n$--dimensional integral submanifold of $\Cc$, or, more 
conceptually, a locally maximal integral submanifold of $\Cc$. ``Locally maximal" means that
even locally $L$ does not belong to an integral submanifold of greater dimension.

These considerations lead to conjecture that 
\begin{quote}
{\it locally maximal integral submanifolds of the Cartan distribution $\Cc^k$ are analogues of
Legendrian submanifolds in contact geometry and that the solutions of a PDE 
$\Ee\subset J^k(E,n)$ are such submanifolds belonging to $\Ee$.}
\end{quote} 

\subsection{Locally maximal integral submanifolds of $\Cc^k$.} 
Motivated by this conjecture we 
shall describe locally maximal integral submanifolds of $\Cc^k$. Let 
$W\subset J^{k-1}(E,n), k\geq1$, 
be an integral submanifold of $\Cc^{k-1}$ which is transversal to the fibers of the projection 
$\pi_{k.k-1}$. By Proposition \ref{C-list}, (6), $\dim\,W\leq n$. Associate with $W$ the 
submanifold $\Ll(W)\subset J^k(E,n)$:
$$
\Ll(W)=\{\theta\in J^k(E,n)\,|\,R_{\theta}\supset T_{\theta'}W \;\mbox{with}\;\theta'=
\pi_{k.k-1}(\theta)\in W\}.
$$
Obviously, $\Ll(L_{(k-1)})=L_{(k)}$ and $\Ll(\{\theta\})=\pi_{k,k-1}^{-1}(\theta)$ for 
any point $\theta\in J^{(k-1)}(E,n)$.
\begin{proc}\label{int-maximal}$($see \cite{V1,KLV,Klava}$)$
\begin{enumerate}
\item $\Ll(W)$ is a locally maximal integral submanifold of $\Cc^k$.
\item If $\dim\,W=s$, then 
$$
\dim\Ll(W)=s+m\left(\begin{array}{c}n+k-s-1 \\n-s-1\end{array}\right).
$$
\item If $N\subset J^k(E,n)$ is a locally maximal integral submanifold, then there is an open 
and everywhere dense subset $N_0$ in $N$ such that 
$$
N_0=\bigcup_{\alpha}U_{\alpha} \quad \mbox{with} \quad U_{\alpha} 
\quad \mbox{being an open domain in}\quad \Ll(W_{\alpha})
$$
\item If $\dim W_1<\dim W_2$, then $\dim \Ll(W_1)>\dim \Ll(W_2)$ except inthe cases $(i)\;
n=m=1, (ii) \;k=m=1$ and $(iii) \;m=1, \,\dim W_1+1=\dim W_2=n$.
\end{enumerate}
\end{proc}
An important consequence of Proposition \ref{int-maximal} is that it disproves the above conjecture.
Existence of locally maximal integral submanifolds  of different dimensions is what makes a  substantial difference between higher order contact structures and  the classical original. In 
particular, this creates a problem in definition of solutions of PDEs in an {\it intrinsic} manner. 
To resolve it we need some additional arguments.

Situations (i)-(iii) in assertion (4) of Proposition \ref{int-maximal} will be called {\it exceptional},
while the remaining ones {\it regular}. This assertion shows that in the regular case integral submanifolds $W_{\alpha}$ figuring in assertion (3) must have the same dimension. This 
dimension will be called the {\it type} of the  maximal integral submanifold $N$. For some other
reasons, which we shall skip, the notion of type can be defined also in exceptional cases (ii) and 
(iii). On the contrary, in the case (i) (classical contact geometry!) all maximal integral submanifolds
are Legendrian and hence are locally equivalent.

Now we may notice that, except for the case $k=m=1$ (classical contact geometry), the fibers 
of the projection $\pi_{k,k-1}, k>1$ are intrinsically characterized as locally maximal integral submanifolds of zero-th type . Therefore, the manifold $J^{k-1}(E,n)$ may be interpreted as
the variety of such submanifolds and, similarly, the distribution $\Cc^{k-1}$  can be recovered
from $\Cc^k$. So, the obvious induction arguments show that by starting from the $k$-th
order contact manifold $(J^{k}(E,n),\Cc^k)$ we can intrinsically recover the whole tower
$$
J^k(E,n)\stackrel{\pi_{k,k-1}}{\longrightarrow} J^{k-1}(E,n)\stackrel{\pi_{k-1,k-2}}
{\longrightarrow}\dots\stackrel{\pi_{\epsilon+1,\epsilon}}{\longrightarrow}J^{\epsilon}(E,n)
$$
where $\epsilon=0$ if $m>1$ and  $\epsilon=1$ if $m=1$. In particular, the projections $\pi_{k,0}$
(resp., $\pi_{k,1}$) can be intrinsically characterized in terms of the $k$--th order contact structure
if $m>1$ (resp., if $m=1$ and $k>1)$. So, if $m>1$, submanifolds $L_{(k)}$ are characterized in these 
terms as locally maximal integral submanifolds of type $n$ that diffeomorphically project on
their images via $\pi_{k,0}$. If $m=1$, then only contact manifold $(J^1(E,n), \Cc^1)$ can be
intrinsically described in terms of a $k$--th order contact structure as the image of the intrinsically
defined projection $\pi_{k,1}$. So, in this case in order to characterize the submanifolds $L_{(k)}$ 
we additionally need to supply the image of $\pi_{k,1}$ with a fiber structure, which mimics  
$\pi_{k,0}$. 

\subsection{Generalized solutions of NPDEs.}
The above considerations lead us to the following definition.
\begin{defi}\label{def-sol}
\begin{enumerate}
\item A locally maximal integral submanifold of type $n$ will be called \emph{$R$-manifold}. 
In particular, submanifolds $L_{(k)}$ are $R$-manifods.
\item Generalized $($resp., ``usual"\,$)$ solutions of a PDE $\Ee\subset J^k(E,n)$ are 
$R$-manifolds $($resp., manifolds $L_{(k)})$ belonging to $\Ee$.
\end{enumerate}
\end{defi}
With this definition we gain 
\begin{quote} {\it the concept of {\bf generalized solutions for nonlinear PDE's}, which, principally, 
cannot be formulated in terms of functional analysis as in the case of linear PDEs.}
\end{quote}
(see \cite{Sob,Schw, GS}). This is one of many instances where a geometrical approach to PDEs
can be in no way substituted by methods of functional analysis or by other analytical  methods.

Definition \ref{def-sol} may be viewed as an extension of the concept of a generalized solution of a linear PDE in the sense of Sobolev-Schwartz to general NPDEs. We have no sufficient 
``space-time" to discuss this very interesting question here. A very rough idea about this relation
is that a generalized solution in the sense of Definition \ref{def-sol} may be viewed as a 
{\it multivalued} one. If the equation is linear, then it is possible to construct a 1-valued one just
by summing up various branches of a multivalued one. The result of this summation is, 
generally, no longer a smooth function but a ``generalized" one. A rigorous formalization of
this idea requires, of course, a more delicate procedure of summation and the Maslov index
(see \cite{Mas}) naturally appears in this context.

\subsection{PDEs versus differential systems.}
According to E.~Cartan, a PDE  $\Ee\subset J^k(E,n)$ can be converted into a differential 
system by restricting the distribution $\Cc^k$ to $\Ee$. The restricted distribution denoted 
by $\Cc^k_{\Ee}$ is defined as
$$
\Cc^k_{\Ee}\,:\,\Ee\ni \theta \mapsto \Cc^k\cap T_{\theta}\Ee.
$$
Originally, E.~Cartan used the Pfaff (exterior) system 
$\omega^i_{\sigma}=0, i=1,\dots,m, |\sigma|<k,$ in order to describe $\Cc^k_{\Ee}$, and this
explains the term {\it exterior differential system}. 

The passage from the equation $\Ee$ understood as a submanifold of $J^k(E,n)$ 
to the differential system  $(\Ee,\Cc^k_{\Ee})$ means, in essence, that we forget that $\Ee$
is a submanifold of $J^k(E,n)$ and consider it as an abstract manifolds equipped with a 
distribution. Cartan was motivated by the idea of replacing {\it non-invariant}, i.e., depending 
on the choice of local coordinate, language of partial derivatives by the invariant calculus of differentials and hence of differential forms. The idea that the general theory of PDEs requires
an invariant and adequate language is of fundamental importance, and E.~Cartan was probably 
the first who raised it explicitly. On the other hand, it turned out later that the language of 
differential forms is not sufficient in this sense. For instance, Proposition \ref{int-maximal} 
 illustrates the fact that the concept of a solution for a generic differential system is not 
well-defined because of the existence of integral submanifolds of different types. The 
{\it rigidity theory} (see \cite{V1,KLV,Klava}) sketched below makes this point  more precise.

First, note that locally maximal integral submanifolds of the restricted distribution $\Cc^k_{\Ee}$
are intersections of such submanifolds for $\Cc^k$ with $\Ee$. So, if $\Ee$ is not very 
overdetermined, i.e., if the codimension of  $\Ee$ in $J^k(E,n)$ is not too big, then the 
difference between locally maximal integral submanifolds of $\Cc^k$ of different types 
survives the restriction to $\Ee$. So, the information about this difference in an explicit 
form gets lost when passing to the differential system $(\Ee,\Cc^k_{\Ee})$. The problem 
to recover it becomes rather difficult especially if the $1$--forms $\omega_i\in\Lambda^1(\Ee)$
of the  Pfaff system $\omega_i=0$ describing the distribution $\Cc^k_{\Ee}$ are arbitrary, 
say, not Cartan ones. Moreover, if we have a generic differential system $(M,\mathcal{D})$ with 
$\mathcal{D}=\{\rho_i=0\}, \,\rho_i\in\Lambda^1(M)$, then it is not even clear which class of its
integral submanifolds should be called solutions. To avoid this inconvenience, E.~Cartan proposed 
to formulate the problem associated with a differential system as the problem of finding its 
integral  submanifolds (locally maximal or not) of a prescribed dimension. But numerous
examples show that a differential system may possess integral submanifolds of an absolutely
different nature, which have the same dimension. One of the simplest examples of this kind
is the differential system $(J^k(E,1),\Cc^k)$ with $\dim E=2, \,k>1$, for which locally maximal
integral submanifolds of types $0$ (fibers or the projection $\pi_{k,k-1}$) and $1$ (R-manifolds)
are all $1$-dimensional. Moreover, integral submanifolds of type $0$ are irrelevant/``parasitic"
in the context of the theory of differential equations.

Secondly, an equation $\Ee\subset J^k(E,n)$ is called {\it rigid} if the $k$-th order contact
manifold $(J^k(E,n,\Cc^k)$ can be recovered if $(\Ee,\Cc^k_{\Ee})$ as an abstract differential
system is only known. For instance, if the codimension of $\Ee$ in $J^k(E,n)$ is less than the difference of dimensions of locally maximal integral submanifolds of types $0$ and $1$, then
$\Ee$ is, as a rule, rigid. Indeed, in this case integral submanifolds of $(\Ee,\Cc^k_{\Ee})$ of 
absolutely maximal dimension are intersections of fibers of $\pi_{k,k-1}$ with
$\Ee$. In other words, these are fibers of the projection 
$\pi_{k,k-1}\left |\right._{\Ee}: \Ee\rw J^{k-1}(E,n)$. If, additionally, this projection is surjective,
then $J^{k-1}(E,n)$ is recovered as the variety of integral submanifolds of 
$(\Ee,\Cc^k_{\Ee})$ of maximal dimension. Next, under some week condition projections of
spaces $\Cc^k_{\Ee,\theta}, \,\theta\in \Ee$, on the so-interpreted jet space $J^{k-1}(E,n)$
span the distribution $\Cc^{k-1}$. In this way  $(J^{k-1}(E,n),\Cc^{k-1})$ is recovered from
$(\Ee,\Cc^k_{\Ee})$ and, finally, $(J^{k}(E,n),\Cc^{k})$ is recovered from $(J^{k-1}(E,n),\Cc^{k-1})$
as the variety of $R$--planes on $J^{k-1}(E,n)$ according to Proposition \ref{C-list}, (1).
Thus converting rigid equations into differential systems is counterproductive, since this 
procedure create non-necessary additional problems.  In this connection it is worth mentioning
that the most important PDEs in geometry, mechanics and physics we deal with are determined or slightly overdetermined systems of PDE's, like Maxwell or Einstein equations, and hence are rigid.
 
Even more important arguments, which do not speak in favor of differential systems, come from 
the fact that the calculus of differential forms is a small part of a much richer structure formed by
{\it natural functors of differential calculus} and {\it  objects representing them}. For instance,
 indispensable for formal integrability theory diff- and jet-Spencer complexes are examples of
 this kind (see \cite{Spen,V-1,V4,Pom,KLV,Klava}). 
Finally, our distrust of differential systems is supported by the fact that practical computations 
of symmetries, conservation laws and other quantities characterizing PDEs become much more complicated in terms of differential systems. 

\subsection{Singularities of generalized solutions.}
The concept of generalized solutions for NPDEs, which is important in itself, naturally leads 
to an important part of a general theory of PDEs, namely, the theory of singularities of 
generalized solutions. Below we shall outline some key points of this theory.

Let $N\subset J^{k}(E,n)$ be an $R$--manifold. A point $\theta\in N$ is called {\it singular of 
type $s$} if the kernel of the differential $d_{\theta}\pi_{k,k-1}$ restricted to $T_{\theta}N$ is of
dimension $s>0$. Otherwise, $\theta$ is called {\it regular}. It should be stressed here that 
``singular" refers to singularities of the map $\pi_{k,k-1}\left |\right._N$, while, by definition, 
$N$ is a smooth submanifold. According to Proposition \ref{C-list}, (6), $N$ is of the form 
$L_{(k)}$ in a neighborhood  of any regular point.

Put $F_{\theta}=(\pi_{k,k-1})^{-1}(\pi_{k,k-1}(\theta))$ (the fiber of $\pi_{k,k-1}$ passing 
through $\theta$) and $V\Cc^k_{\theta}=\Cc^k_{\theta}\cap T_{\theta}\left(F_{\theta}\right)$.  
The {\it bend} of $N$ at a point $\theta\in N$ is
$$
B_{\theta}N\df\ker d_{\theta}\left(\pi_{k,k-1}\left |\right._N\right)=
T_{\theta}N\cap T_{\theta}\left(F_{\theta}\right)\subset V\Cc^k_{\theta}.
$$
Also, we shall call an {\it $s$--bend} (at $\theta\in J^k(E,n)$) an $s$--dimensional subspace 
of $V\Cc^k_{\theta}$, which is of the form $B_{\theta}N$ for some $R$--manifold $N$. Bends 
are very special subspaces in $V\Cc^k_{\theta}$. A remarkable fact is that $s$--dimensional 
bends are classified by $s$--dimensional Jordan algebras of a certain class over $\dR$, which contains all unitary algebras (see \cite{Vsing,V<}).
\begin{quote}
{\it PDEs differ from each other by the types of singu-\\larities   which their generalized 
solutions admit.}
\end{quote}
For instance, 2-dimensional Jordan algebras associated with 2-dimensional bends are 
2-dimensional unitary algebras and hence are isomorphic to one of the following three algebras
$$
\dC_{\epsilon}=\{a+b\zeta\,|\,a,b\in\dR, \,\zeta^2=\epsilon 1\} \quad 
\mbox{with} \quad \epsilon=\pm \;\mbox{or} \;0.
$$
Obviously, $\dC_{-}=\dC$ and $\dC_{+}=\dR\oplus\dR$ (as algebras). An equation in two 
independent variables  is {\it elliptic} (resp., {\it parabolic} or {\it hyperbolic}) if its generalized
solutions possess singularities of type $\dC_{-}$ (resp., $\dC_{0}, \,\dC_{+}$) only. Geometrically, 
singularities corresponding to algebra $\dC$ are Riemann ramifications, while bicharacteristics
of hyperbolic equations reflect the fact that $\dC_{+}$ splits into the direct sum $\dR\oplus\dR$.

Obviously, the simplest singularities correspond to the algebra $\dR$. They present a kind of 
folding and can be {\it analytically detected} in terms of non-uniquiness of Cauchy data. 
A similar analytic approach is 
hardly possible for more complicated algebras. This explains why analogues of the classical subdivision of PDEs in two independent variables into elliptic, parabolic and hyperbolic ones 
are not yet known. This fact emphasizes once again that only analytical methods for PDEs, 
even linear ones, are not sufficient  and the geometrical approach is indispensable.

Description of singularities that solutions of a given PDE admit is naturally settled as follows. 
Let $\Sigma$ be a type of $s$--bends, which may be identified with the corresponding
Jordan algebra. If $N$ is an $R$--manifold, then
$$
N_{\Sigma}=\{\theta\in N\,|\,B_{\theta}N\;\mbox{is of the type}\;\Sigma\}
$$
is the locus of its singular points of type $\Sigma$. Generally, $\dim N_{\Sigma}=n-s$. If $N$
is a solution of a PDE $\Ee$, then $N_{\Sigma}$ must satisfy an auxiliary system of PDEs,
which we denote by $\Ee_{\Sigma}$.  For "good" equations $\Ee_{\Sigma}$ is, generally, 
a nonlinear, undetermined system of PDEs in $n-s$ independent variables.

\subsection{The reconstruction problem.}
So, any PDE is not a single but is surrounded by an ``aura" of subsidiary equations,
which put in evidence the {internal structure} of its solutions. The importance of these equations 
becomes especially clear in the light of the {\it reconstruction problem}:
\begin{quote}
{\it Whether the behavior of singularities of solutions of a PDE $\Ee$ \\ uniquely determines the equation itself or, equivalently, whether \\ is possible to reconstruct  $\Ee$ assuming that the
$\Ee_{\Sigma}$'s are known?} 
\end{quote}
In a physical context this question sounds as
\begin{quote}
{\it Whether the behavior of singularities of a field (medium, etc.) \\ completely determines 
the behavior of the field (medium, etc.) itself?}
\end{quote}
A remarkable example of this kind is the deduction of Maxwell's equations from elementary
laws of electricity and magnetism (Coulomb, \dots, Faraday) (see \cite{Lun}).

The reconstruction problem resolves positively for hyperbolic NPDEs on the basis of 
equations $\Ee_{\mathrm{FOLD}}$ that describes singularities corresponding to the algebra 
$\dR$.  The equations describing wave fronts of solutions of a linear hyperbolic PDE $\Ee$ 
are part of the system $\Ee_{\mathrm{FOLD}}$.
\begin{ex}\label{mu}
\emph{Fold--type singularities for the equation $u_{xx}-\frac{1}{c^2}u_{tt}-mu^2=0$.}\par
Consider wave fronts of the form $x=\varphi(t)$ and put
$$
g=u|_{\textrm{wave front}},\quad h=u_x|_{\textrm{wave front}}.
$$
Then we have
$$\left\{\begin{array}{c}\ddot{g}+(cm)^2g=\pm 2c\dot{h} \\1-\frac{1}{c^2}\dot{\varphi}^2=0\Leftrightarrow\dot{\varphi}=\pm c\end{array}\right.\Leftarrow\left[\begin{array}{c}\textrm{Equations describing}  \;the \\\textrm{behavior of fold--type} \\\textrm{singularities}\end{array}\right.$$
The second of these equations is of eikonal type and describe the space-time shapes of 
singularities. On the contrary, the first equation describes a ``particle" in the ``field" $h$.
If this field is constant $\Leftrightarrow \;\dot{h}=0$, then the first equation represents a 
harmonic oscillator of frequency $\nu=mc$.
\end{ex}
\begin{ex}\label{Kl-Gord}
\emph{Fold--type singularities for the Klein--Gordon equation 
$$(\partial_t^2 -\vec{\nabla}^2+m^2)u=0.$$}
Consider wave fronts of the form $t=\varphi(x_1,x_2,x_3)$ and $g$ and $h$ as in example \ref{mu}\\
$\Ee_{\emph{FOLD}}$ = $\left\{\begin{array}{l}(\vec{\nabla}\varphi)^2=1\leftarrow\quad\textrm{eikonal type equation} \\ \nabla^2h+m^2h-g-(\nabla^2\varphi)g=2\vec{\nabla}\varphi\cdot\vec{\nabla}g\leftarrow\textrm{???}\end{array}\right.$\par
The physical meaning of the second of these equations is unclear.
\end{ex}
\begin{ex}
Classical Monge-Amp\`{e}re equations are defined as equations of the form
$$
S(u_{xx}u_{yy} - u_{xy}^2)+ Au_{xx} + Bu_{xy} + Cu_{yy} + D = 0 
$$
with $S, A, B, C, D$ being functions of $x, y, u, u_x, u_y$ $($see \cite{KLR}$)$. As it was already observed by S.~Lie this class of equations is invariant with respect to contact transformations. 
This fact forces to think that Monge-Amp\`ere equations are distinguished by some ``internal"
property.  This is the case, and  Monge-Amp\`ere equations are completely characterized by
the fact that the reconstruction problem for these equations is equivalent
to a problem in contact geometry $($see \cite{CFV,MVY}$)$. 
\end{ex}
 The reader will find in \cite{Eik} further details and examples concerning the auxiliary singularities equations. Some exact generalized solutions  of Einstein equations (the ``square root" of the 
Schwarzshild solution, etc) are described in \cite{SVV}.  

\subsection{Quantization as a reconstruction problem.} 
Let $\Ee$ be a PDE, whose solutions admit fold-type singularities. Then we have the 
following series of interconnected equations:
\begin{equation}\label{char}
\Ee \Longrightarrow \Ee_{\mathrm{FOLD}}\Longrightarrow 
\Ee_{\mathrm{eikonal}} \Longrightarrow \Ee_{\mathrm{char}}.
\end{equation}
Here $\Ee_{\mathrm{eikonal}}$ is the equations from the system $\Ee_{\mathrm{FOLD}}$ 
that describes space-time shape (``wave front") of fold-type singularities. It is a Hamilton-Jacobi
equation (see Examples \ref{mu} and \ref{Kl-Gord}). In its turn $\Ee_{\mathrm{char}}$
is the system of ODEs that describes {\it characteristics} of $\Ee_{\mathrm{eikonal}}$. In the 
context where space-time coordinates are independent variables $\Ee_{\mathrm{eikonal}}$ 
is a Hamiltonian system whose Hamiltonian is the main symbol of $\Ee$. Now
we see that the correspondence
\begin{equation}\label{CHAR}
\mbox{CHAR}\,: \Ee \,(\mathrm{PDE}) \Longrightarrow \Ee_{\mathrm{char}} 
\,(\mathrm{Hamiltonian \;system \;of \;ODE's})
\end{equation}
is parallel to the correspondence between quantum  and classical mechanics
\begin{equation}\label{BOHR}
\mbox{BOHR}\,: (\mathrm{Schroedinger's \;PDE}) \Longrightarrow 
\,(\mathrm{Hamiltonian \;ODEs}).
\end{equation}
Moreover, the correspondence (\ref{CHAR}) is at the root of the famous ``optics-mechanics 
analogy", which guided E.~Schr\"odinger in his discovery of the ``Schr\"odinger equation"
(see Schr\"odinger's Nobel lecture \cite{Schr}).

It is remarkable that in ``Cauchy data" terms correspondence  (\ref{CHAR})  was known already 
to T.~Levi-Civita and he tried to put it at the foundations of quantum mechanics (see \cite{L-C}).
From what is known today this attempt was doomed to failure. However, the idea that 
quantization is something like the reconstruction problem explains well why numerous 
quantization procedures proposed up to now form a kind of recipe
book not based on some universal principles. Indeed, from this point of view the quantization
looks like an attempt to restore the whole system $\Ee_{\mathrm{FOLD}}$ on the basis of
knowledge of $\Ee_{\mathrm{char}}$ only. This is manifestly impossible, since 
$\Ee_{\mathrm{char}}$ depends only on the main symbol of $\Ee$. On the other hand, the 
above outlined solution singularity theory admits some interesting generalizations and refinements, which not only keep alive the Levi-Civita idea but even make it more attractive.

\subsection{Higher order contact transformations and the Erlangen program.} 
The above interpretation 
of PDEs as submanifolds of higher order contact manifolds is the first step toward a ``conceptualization" of the standard approach to PDE's. It is time now to test its validity 
through the philosophy of the Erlangen program. First of all, this means that we have to 
describe the symmetrie group of higher contact geometries, i.e., the group of {\it higher 
contact transformations}.
\begin{defi}
A diffeomorphism/transformation $\Phi:J^k(E,n)\rw J^k(E,n)$ is called a \emph{$k$-order contact} 
if for any $X\in \Cc^k, \,\Phi(X)\in\Cc^k$ or, equivalently, 
$d_{\theta}\Phi(\Cc^k_{\theta})=\Cc^k_{\Phi(\theta)}, \,\forall \theta\in \Phi$.
\end{defi}
If $\Phi$ is a $k$--th order contact, then, obviously, it preserves the class of locally maximal 
integral submanifolds of type $s$. In particular, it preserves fibers of the projection $\pi_{k,k-1}$ 
and hence locally maximal integral submanifolds of type $n$ that are transversal to these fibers. 
But the latter are locally of the form $L_{(k)}$ (Proposition \ref{C-list}, (6)). This proves that the differential of $\Phi$  sends $R$--planes into $R$--planes. By identifying these $R$--planes with points of $J^{k+1}(E,n)$ we see that $\Phi$ induces a diffeomorphism $\Phi_{(1)}$ of $J^{k+1}(E,n)$.
More exactly, if $\theta\in J^{k+1}(E,n)$ and $\theta'=\pi_{k+1,k}(\theta)$, then 
$(d_{\theta'}\Phi)(R_{\theta})$ is an $R$-plane and hence is of the form $R_{\vartheta}$ for
a $\vartheta\in J^{k+1}(E,n)$. Then we put $\Phi(\theta)=\vartheta$. Moreover, it directly follows 
from Proposition \ref{C-list}, (1), that $\Phi_{(1)}$ is a $(k+1)$--order contact and the diagram
$$
\begin{array}{ccc}
J^{k+1}(E,n)&\stackrel{\Phi_{(1)}}{\longrightarrow}&J^{k+1}(E,n) \\
\downarrow\lefteqn{\pi_{k+1,k}} & &\downarrow\lefteqn{\pi_{k+1,k}}\\
J^{k}(E,n)&\stackrel{\Phi}{\longrightarrow}&J^{k}(E,n)
\end{array}
$$
commutes. By continuing this process we, step by step, construct contact transformations
$$
\Phi_{(l)}:J^{k+l}(E,n)\stackrel{F_{(1)}}{\longrightarrow}J^{k+l}(E,n), \quad \Phi_{(l)}\df
(\Phi_{(l-1)})_{(1)}.
$$
\begin{thm}\label{L-B}
Let  $\Phi\,:\,J^{k}(E,n)\to J^{k}(E,n), \,k>0,$ be a $k$--order contact transformation.
Then $\Phi=\Psi_{(l)}$ $($resp., $\Phi=\Psi_{(l-1)}$ $)$ where $\Psi$ is a diffeomorphism 
of $E$ if $m>1$ $($resp., a contact transformation of $J^{1}(E,n)$ if $m=1$ $)$.
\end{thm}
A proof of this fundamental result for the classical symmetry theory can be easily
deduced from the fact explained above that a $k$--th order contact transformation 
preserves fibers of $\pi_{k,k-1}$ and hence induces a $(k-1)$--th order contact 
transformation of $J^{k-l}(E,n)$. For $m=1$ it was proven by Lie and B\"acklund 
(see \cite{V1,KVin}).

If one takes Definition  \ref{eq.jet} for a true definition of PDEs, then the definition of a 
symmetry of a PDE should be
\begin{defi}\label{sm}
A symmetry of a PDE $\Ee\subset J^k(E,n)$ is \\
\emph{(1)} a $k$--th order contact transformation $\Phi: J^k(E,n)\to J^k(E,n)$ such that 
$\Phi(\Ee)=\Ee$  $($\`a la S.~Lie$)$;\\
\emph{(2)} a diffeomorphism $\Psi:\Ee\to \Ee$ preserving the distribution $\Cc_{\Ee}^k$
$($\`a la E.~Cartan$)$.
\end{defi}
The rigidity theory shows that Definitions (1) and (2) are equivalent for rigid PDEs, i.e., for 
almost all PDEs of practical interest. Moreover, by Theorem \ref{L-B}, Definitions \ref{sm} 
and \ref{cl-sym} are equivalent in this case too.

\begin{rmk}There are analogues of theorem \ref{L-B} and definition \ref{sm} for infinitesimal 
$k$--order contact transformations and symmetries. They do not add anything new 
to our discussion, and we shall skip them.
\end{rmk}

In the light of the ``Erlangen philosophy" the result of Theorem \ref{L-B} looks disappointing.
Indeed, it tells us that the group of $k$-order contact transformations coincides with the group
of first order transformations. So, higher order contact geometries are governed by the
same group as the classical one. This does not meet a natural expectation that 
transformations of higher order geometries should form some larger groups. Hence, by 
giving credit to this philosophy, we are forced to conclude that
\begin{quote}
{\it  Definition \ref{eq.jet} or what is commonly meant by a differential equation is not 
a conceptual definition but should be considered just as a description of an object, 
whose nature must be still discovered.}
\end{quote}
So, the question of what object is hidden under this description is to be investigated.
One rather evident hint is to examine the remaining case $k=\infty$. This is psychologically
difficult, since $J^{\infty}(E,n)$ being an infinite-dimensional manifold of a certain kind does 
not  possess any ``good" topology or norm, etc. which seem indispensable for the existence of a 
``good" differential calculus on it. Another hint comes from  the principle ``chercher la sym\'etrie".  
For instance, if $\Ee$ (resp., $\square$) is a linear equation (resp., a linear differential operator) 
with constant coefficients, then  $\square$ sends solutions of $\Ee$ to the solutions.  
For this reason $\square$ may be considered as a symmetry of $\Ee$, finite or infinitesimal. Symmetries of this kind are not, generally, classical and their analytical description involves 
partial derivatives of any order. Hence one may expect that something similar takes place for
general PDEs, and we are going to show that this is the case.

\section{From integrable systems to diffieties and higher symmetries.}
\subsection{New experimental data: integrable systems.}
The discovery in the late 1960s of some remarkable properties of the now famous 
Korteweg - de Vries equation and later of other {\it integrable systems} brought to light 
various new facts, which had no conceptual explanation in terms of the classical symmetry
theory. In particular, any such equation is included in an infinite series of similar equations,
the {\it hierarchy}, which are interpreted as commuting Hamiltonian flows with respect to an, in 
a sense, infinite-dimensional Poisson structure. For this reason equations of this hierarchy
may be considered as infinitesimal symmetries of each other. Moreover, they involves derivatives
of any order and hence are outside the classical theory (see \cite{ZF}).  So, attempts to include
these non-classical symmetries into common with classical symmetries frames directly leads
to infinite jets.

\subsection{Infinite jets and infinite order contact transformations.} Recall that the Cartan distribution
$\Cc^{\infty}$ on $J^{\infty}(E,n)$ is (paradoxically!) $n$--dimensional and completely integrable
(Proposition \ref{C-list}, (5)). A consequence of this fact is that locally maximal integral submanifolds 
of $\Cc^{\infty}$ are of the same type in sharp contrast with finite-order contact geometries
(Proposition \ref{int-maximal}). This is a weighty argument in favor of infinite jets. After that we
have to respond to the question of whether the group of infinite-order contact transformations is 
broader than the group of classical ones. More exactly, we ask whether there are infinite-order 
contact transformations that are not of the form $\Phi_{(\infty)}$ where $\Phi$ is a finite-order 
contact transformation (see Theorem \ref{L-B}). Here $\Phi_{(\infty)}$ stands for the direct limit 
of $\Phi_{(l)}$'s. The answer is positive: this (local) group consists of all {\it invertible differential operators} (in the generalized sense outlined above) acting on $n$--dimensional submanifolds of 
$E$. These operators involve partial derivatives of arbitrary orders and in this sense they justify the 
credit given to infinite jets. We shall skip the details (see \cite{Vloc}), since the same question 
about infinite order infinitesimal symmetries is much more interesting from the practical point of 
view and at the same time it reveals some unexpected a priori details, which become essential for 
the further discussion.

Recall that an infinitesimal symmetry of a distribution $\Cc$ on a manifold $M$ is a vector 
field $X\in D(M)$ such that $[X,Y]\in \Cc$ if $Y\in \Cc$ (symbolically, $[X,\Cc]\subset \Cc$). 
Infinitesimal symmetries form a subalgebra in $D(M)$ denoted $D_{\Cc}(M)$. The flow
generated by a field $X\in D_{\Cc}(M)$ moves, if it is globally defined, (maximal) integral 
submanifolds of $\Cc$ into themselves. If it not globally defined this flow moves only sufficiently
small pieces of integral submanifolds. In this sense we can speak of a local flow in the ``space
of  (maximal) integral submanifolds of  $\Cc$\,".

If the distribution $\Cc$ is integrable/Frobenius, then it may be interpreted as a foliation whose 
leaves are its locally maximal integral submanifold. In this case $\Cc$ is an ideal in $D_{\Cc}(M)$. 
If $N\subset M$ is a leaf of  $\Cc$, then any $Y\in \Cc$ is tangent to $N$ and, therefore, the 
flow generated by $Y$ leaves $N$ invariant, i.e., any leaf of  $\Cc$ slides along itself under 
the action of this flow. We may interpret this fact by saying that the local flow generated by $Y$  
on the ``space of all leaves of $\Cc$\," is trivial. This is, obviously, no longer so if 
$Y\in D_{\Cc}(M)\setminus\Cc$.  Hence the flow generated by $Y$ in the ``space of all leaves 
of $\Cc$\," is uniquely defined by the coset $[Y \mod \Cc]$, and the quotient Lie algebra
\begin{equation}\label{q}
\Sym\Cc\df\frac{D_{\Cc}(M)}{\Cc},
\end{equation}
called the {\it symmetry algebra of $\Cc$}, is naturally interpreted as the algebra of {\it vector 
fields on the ``space of leaves\," of} $\Cc$. It should be stressed that it would be rather counterproductive to try to give a rigorous meaning to  the ``space of leaves\,". On the 
contrary, the above interpretation of the quotient algebra (\ref{q}) is  very productive and 
may be interpreted as the smile of the Cheshire Cat.

Now we shall apply the above construction to the distribution $\Cc^{\infty}$ and introduce for 
this special case the following notation:
\begin{equation}\label{C-note}
\Cc D(J^{\infty}(E,n))=\Cc^{\infty},  \quad D_{\Cc}(J^{\infty}(E,n))=D_{\Cc^{\infty}}(J^{\infty}(E,n)), 
\quad \varkappa=\Sym \Cc^{\infty}. 
\end{equation}
The Lie algebra $\varkappa$ will play a prominent role in our subsequent investigation. At the 
moment we know that it is the ``algebra of vector fields on the space of all locally maximal 
integral submanifolds of $J^{\infty}(E,n)$\,". As a first step we have to describe $\varkappa$ in
coordinates. However, in order to do that with a due rigor we have to clarify  
before what is differential calculus on infinite-dimensional manifolds of the kind. It is rather
obvious that the usual approaches based on ``limits", ``norms", etc, cannot be applied to this 
situation. So, we need the following digression.

\subsection{On differential calculus over commutative algebras.}\label{dif-cal} Let $A$ be 
a unitary, i.e., commutative and with unit, algebra over a field $\boldsymbol{k}$ and $P$ 
and $Q$ be some $A$--modules. 
\begin{defi}\label{def-DO}
$\Delta:P\longrightarrow Q$ is a \emph{linear differential operator (DO)} of order $\leq m$ if $\Delta$ is 
$\gk$--linear and $[a_0,[a_1,\ldots,[a_m,\Delta]\ldots]]=0$, $\forall a_0,a_1,\ldots, a_m\in A$.
\end{defi}
\noindent Elements $a_i\in A$ figuring in the above multiple commutator are understood as the multiplication by $a_i$ operators.

If $A=C^{\infty}(M), \,P=\Gamma(\pi), \,Q=\Gamma(\eta)$ with $\pi, \eta$ being some vector 
bundles, then  Definition \ref{def-DO} is equivalent to the standard one. The ``logic" of differential calculus is formed by {\it functors of differential calculus} together with their natural 
transformations and representing them objects in a {\it differentially closed} category of 
$A$--modules \cite{V-1,V1,V4}. In particular, this allows one to construct analogues of all 
known structures in differential geometry, say, tensors, connections, de Rham and Spencer cohomology, an so on, over an arbitrary unitary algebra. The reader will find in \cite{Nestr} an elementary introduction to this subject based on a {\it physical motivation}.

By applying this approach 
to the filtered algebra $\Ff_{\infty}=\{\Ff_i\}$ (see (\ref{F-tower})) we shall get all necessary 
instruments to develop differential calculus on spaces $J^{\infty}(E.n)$ and, more generally, on 
{\it diffieties} (see below). The informal interpretation of the filtered algebra $\Ff_{\infty}$ as the smooth function algebra on the ``cofiltered manifold" $J^{\infty}(E.n)$ helps to keep the analogy 
with the calculus on smooth manifolds under due control. In coordinates an element 
of $\Ff_{\infty}$ looks as a function of a finite number of variables $x_i$ and $u^j_{\sigma}$.
This reflects the fact that any  ``smooth function" on  $J^{\infty}(E.n)$ is, by definition, a
smooth function on a certain $J^k(E.n), \,k<\infty$, pulled back onto $J^{\infty}(E.n)$ via $\pi_{\infty,k}$.
So, the filtered structure of $\Ff_{\infty}$ is essential and differential operators
$\Delta:\Ff_{\infty}\longrightarrow \Ff_{\infty}$ (in the sense of definition \ref{def-DO}) must 
respect it. This means that $\Delta(\Ff_k)\subset \Ff_{k+s}$ for some $s$. In particular, a vector 
field on $J^{\infty}(E.n)$ is defined as a derivation of  $\Ff_{\infty}$, which respects, in this sense,
 the filtration. In coordinates such a vector field looks as an infinite series
\begin{equation}\label{vf-on-jets}
 X=\sum_i\alpha_i\frac{\p}{\p x_i}+\sum_{j, \sigma}\beta_{\sigma}^j\frac{\p}{\p u_{\sigma}^j}, 
 \quad \phi_i,\,\psi_{\sigma}^j\in \Ff_{\infty}.
\end{equation}
The $\Ff_{\infty}$--module of vector fields on $J^{\infty}(E.n)$ will be denoted 
by $D(J^{\infty}(E.n))$.

 \subsection{Algebra  $\vk$ in coordinates.} 
 Since $\Cc^{\infty}$ is an $\Ff_{\infty}$--module generated by the vector fields $D_i$'s 
 (Proposition \ref{C-list}, (5)), it is convenient 
 to represent a vector field $X\in D(J^{\infty}(E.n))$ in the form
 \begin{equation}\label{vf-on-jets}
 X=\sum_i\psi_iD_i+\sum_{j, \sigma}\f_{\sigma}^j\frac{\p}{\p u_{\sigma}^j}, 
 \quad \psi_i,\,\f_{\sigma}^j\in \Ff_{\infty}.
\end{equation}
where the first summation, which belongs to $\Cc^{\infty}$, is the {\it horizontal} part of $X$, while
the second one is its {\it vertical} part. This splitting of a vector field into horizontal and vertical
parts is unique but depends on the choice of  coordinates. Obviously, the coset  
$[X \mod \Cc^{\infty}]$ is uniquely characterized by the vertical part of $X$.

Below we use the notation $D_\sigma\df D_1^{\sigma_1}\cdot\dots\cdot D_n^{\sigma_n}$ for
a multiindex $\sigma=(\sigma_1,\dots,\sigma_n)$.
\begin{proc}
\begin{enumerate}
\item $\vk$ is a $\Ff_{\infty}$--module and $\p/\p u^1,\dots,\p/\p u^m$ is its local basis in the
chart $U$ with coordinates $(\dots,x_i,\dots,u^j_{\sigma},\dots)$;
\item the correspondence
$$
(\Ff_{\infty}^m)_U\ni \varphi=(\varphi^1,\ldots, \varphi_m) \Leftrightarrow \ev_\varphi=\sum D_\sigma(\varphi_i)\frac{\partial}{\partial u_\sigma^i}\in \vk_U
$$
is an isomorphism of $\Ff_{\infty}$--modules localized to the chart $U$;
\item the Lie algebra structure $\{\cdot,\cdot\}$ in $\vk_U$ is given by the formula
$$
\{\varphi,\psi\}=\ev_\varphi(\psi)-\ev_\psi(\varphi), \quad  [\ev_\varphi,\ev_\psi]=\ev_{\{\varphi,\psi\}};
$$
\item $(f,\ev_{\varphi})\mapsto \ev_{f\varphi}$ is the $(\Ff_{\infty})_U$--module product in $\vk_U$.
\end{enumerate}
\end{proc}
The vector fields $\ev_{\f}$s locally representing elements of the module $\vk$ are called 
{\it evolutionary derivations}, and $\f$ is called the {\it generating function} of $\ev_{\f}$. The 
bracket $\{\cdot,\cdot\}$ introduced for the first time in \cite{Vsym} (see also \cite{Vloc,KVin}) 
is a generalization of both the Poisson and the contact brackets. Indeed, these are particular
cases where $m=1$ and the generating functions depend only on the $x_i$s and on the first derivatives and in the contact case also of  $u$. If $Y$ is a vector field on $E$ ($m>1$) or a 
contact vector field on $J^1(E,n)$ ($m=1$) and  $Y_{(\infty)}$ is its lift to $J^{\infty}(E,n)$, 
then $Y_{(\infty)}\in D_{\Cc}(J^{\infty}(E,n))$ and the composition 
$Y\mapsto  Y_{(\infty)}\mapsto [Y_{(\infty)} \mod \Cc D(J^{\infty}(E,n))]\in\vk$ is injective. So,
infinitesimal point and contact transformations are naturally included in $\vk$. Their generating 
functions depends only on $x, u$ and first derivatives, and we see that the Lie algebra $\vk$ is
much larger than the algebras of infinitesimal point and contact transformations. Hence the
passage to infinite jets is in fairly good accordance with the ``Erlangen philosophy". But in order 
to benefit from this richness of infinite order contact transformations we must bring PDEs in the 
context of infinite order contact geometry. 
But in that case we cannot mimic Definition \ref{eq.jet}, since, in sharp contrast with
finite order jet spaces, an arbitrary submanifold of $S\subset J^{\infty}(E,n)$ can not be interpreted 
as a PDE. Indeed, the restriction of $\Cc^{\infty}$ to $S$ is, generally, not $n$--dimensional, while
we need $n$--dimensional integral submanifolds to define the solutions. So, we must concentrate
on those submanifolds $S$ to which $\Cc^{\infty}$ is tangent, i.e., such that
$\Cc^{\infty}_{\theta}\subset T_{\theta}S, \,\forall \theta\in S$. These are obtained
by means of the {\it prolongation procedure}.

\subsection{Prolongations of PDEs and diffieties.} 
Let $\Ee\subset J^k(E,n)$ be a PDE in the sense of Definition 
\ref{eq.jet} and $N\subset \Ee$ be its solution (Definition \ref{def-sol}). Then, obviously, 
$T_{\theta}N\subset \Ee_{\theta}, \,\forall _\theta\in N$. So, if $\Ee$ admits a solution passing
through a point $\theta\in\Ee$, then there is at least one $R$--plane at $\theta$, which is tangent
to $\Ee$. Since any $R$--plane is of the form $R_{\vartheta},  \pi_{k+1,k}(\vartheta)=\theta$, the
variety of all $R$--planes tangent to $\Ee$ is identified with the submanifold (probably, with singularities)
$$
\Ee_{(1)}\df\{\vartheta\in J^{k+1}{(E,n)}\,|\,R_{\vartheta} \;\mbox{is tangent to} \;\Ee\}
\subset J^{k+1}(E,n).
$$
So, tautologically, a solution of $\Ee$ passes only through points of 
$\pi_{k+1,k}(\Ee_{(1)})\subset \Ee$. In other words, a solution of $\Ee$ is automatically a 
solution of $\pi_{k+1,k}(\Ee_{(1)})$. Hence
by substituting  $\pi_{k+1,k}(\Ee_{(1)})$ for $\Ee$ we eliminate ``parasitic" points. Moreover,
by construction, if $L_{(k)}\subset\Ee$, then $L_{(k+1)}\subset\Ee_{(1)}$ and vice versa. 
Hence $\Ee$ and $\Ee_{(1)}$ have common ``usual" solutions but $\Ee_{(1)}$ is without
``parasitic" points of $\Ee$. By continuing this process of elimination of ``parasitic" points we
inductively construct successive {\it prolongations} $\Ee_{(r)}\df\left(\Ee_{(r-1)}\right)_{(1)}$ 
of $\Ee$. In this way we get an infinite series of equations, which have common ``usual" solutions:
\begin{equation}\label{prolong}
\Ee=\Ee_{(0)}\stackrel{\pi_{k+1,k}}{\longleftarrow}\Ee_{(1)}\stackrel{\pi_{k+2,k+1}}{\longleftarrow}\Ee_{(2)}\stackrel{\pi_{k+3,k+2}}{\longleftarrow}\ldots, 
\quad\mbox{with} \quad \Ee_{(r)}\subset J^{k+r}.
\end{equation}
The inverse limit $\Ee_{\infty}$ of the sequence (\ref{prolong})  called the {\it infinite prolongation}
of $\Ee$ is a submanifold of $J^{\infty}(E,n)$
(in the same sense as the latter) and one of the results of the {\it formal theory} of PDE's tells:
\begin{proc}
If the distribution $\Cc^{\infty}$ is tangent to a submanifold $S\subset J^{\infty}(E,n)$, then
$S=\Ee_{\infty}$ for a PDE $\Ee$. 
\end{proc}
\noindent(see \cite{Jan,Spen,Pom,KLV}).

In coordinates, prolongations of $\Ee$ are described as follows
\begin{equation}\label{prolong-coord}
\Ee_{(2)}=\left\{\begin{array}{c}\Ee_{(1)}=\left\{\begin{array}{c}\Ee=\{ F_s(x,u,\ldots,u_\sigma^j,\ldots)=0,\ s=1,\ldots,l\} \\D_i F_s=0\end{array}\right\} \\D_iD_j F_s=0\end{array}\right\}
\end{equation}
\centerline{$\dots$}
\centerline{$\Downarrow$}
\centerline{$\Ee_{\infty}=\{D_\sigma F_s=0, \forall s,\sigma\}$}
\begin{rmk}
$\Ee_{\infty}$ may be empty.
\end{rmk}
The algebra $\Ff_{\infty}(\Ee)\df\Ff|_{\Ee_{\infty}}$ plays the role of the smooth function 
algebra on $\Ee_{\infty}$. It is a filtered algebra
\begin{equation}\label{E-funct} 
\Ff_0(\Ee)\subset\cdots\subset\Ff_s(\Ee)\subset\cdots\Ff_{\infty}(\Ee)
\; \mbox{with}\, \;\Ff_s(\Ee)=\im(C^{\infty}(\Ee_{(s)})
\stackrel{\pi_{\infty,k+s}^*}{\longrightarrow}\Ff_{\infty}(\Ee)).
\end{equation}
As in the case of infinite jets differential calculus on $\Ee_{\infty}$ is understood as
differential calculus over the filtered algebra $\Ff_{\infty}(\Ee)$. 

Thus we have constructed the central object of general theory of PDEs.
\begin{defi}
The pair $(\Ee_{\infty}, \Cc^{\infty}_{\Ee})$ with $\Cc^{\infty}_{\Ee}\df
\Cc^{\infty}\left |\right._{\Ee_{\infty}}$ is called the \emph{diffiety} associated with $\Ee$.
\end{defi}
The distribution $\Cc^{\infty}_{\Ee}$ is $n$--dimensional, since $\Cc^{\infty}$ is tangent to 
$\Ee_{\infty}$. The projection $\pi_{\infty,k}$ establishes a one-to-one correspondence between
integral submanifolds of $\Cc^{\infty}_{\Ee}$ and those of $\Cc^k\left |\right._{\Ee}$, which are
transversal to fibers of $\pi_{k,k-1}$. So, $n$--dimensional integral submanifolds of 
$\Cc^{\infty}_{\Ee}$ are identified with non-singular solutions of $\Ee$. 

The following interpretation, even though absolutely informal, is a very good guide in the task
of deciphering the native language that NPDEs speak and, therefore, in terms of which they can 
be only understood adequately:
\begin{quote}
\it {The diffiety associated with a PDE $\Ee$ (in the standard sense of this term) is the 
\emph{space of all solutions of  $\Ee$}}.
\end{quote}
\begin{rmk}
The reader may have already observed that nontrivial generalized solutions of $\Ee$ 
cannot be interpreted as integral submanifolds of $\Cc^{\infty}_{\Ee}$ and hence the diffiety
$(\Ee,\Cc^{\infty}_{\Ee})$ is not the ``space of all solutions of  $\Ee$". However, this is not 
a conceptual defect, since this diffiety can be suitably completed.
\end{rmk}
As a rule, diffieties are infinite-dimensional.  Diffieties of finite dimension are foliations, probably,
with singularities. Diffieties associated with determined and overdetermined systems of ordinary differential equations (ODEs) are 1--dimensional foliations on finite-dimensional manifolds. On 
the contrary, diffieties associated with underdetermined systems of ODEs are infinite-dimensional.
A good part of control theory is naturally interpreted as structural theory of this kind of diffieties
(see \cite{Fls}).

\subsection{Higher infinitesimal symmetries of PDEs.} \label{High-sym}
Now having in hands the concept of  diffiety 
we can extend the classical symmetry theory described above  by including in it the new already 
mentioned ``experimental data" that come from the theory of integrable systems. To this
end it is  sufficient to apply the same approach we have used to understand what are 
infinite-order infinitesimal contact transformations. 

As before, by abusing the language, we shall denote the $\Ff_{\infty}(\Ee)$--module of vector 
fields on $\Ee_{\infty}$ belonging to $\Cc^{\infty}_{\Ee}$ by the same symbol $\Cc^{\infty}_{\Ee}$.
Since the distribution $\Cc^{\infty}$ is tangent to $\Ee_{\infty}$, vector fields $D_i$s are also
tangent to $\Ee_{\infty}$. For this reason restrictions of the $D_i$ to $\Ee_{\infty}$ are well-defined,
Denote them by $\bar{D}_i$. The Lie algebra of infinitesimal transformations preserving the
distribution $\Cc^{\infty}_{\Ee}$ is
\begin{equation}\label{DC}
D_{\Cc}(\Ee_{\infty})\df\{X\in D(\Ee_{\infty})\,|\,[X,Y]\in \Cc^{\infty}_{\Ee}, 
\,\forall Y\in \Cc^{\infty}_{\Ee}\}.
\end{equation}
Now the Lie algebra of {\it infinitesimal higher symmetries} of a PDE $\Ee$ is defined as
\begin{equation}\label{sym}
\boxed{\Sym\Ee=\frac{D_{\Cc}(\Ee_{\infty})}{\Cc_{\Ee}^{\infty}}}
\end{equation}
This definition merits some comments. First, we use the adjective ``higher" to stress the fact 
that generating functions of elements of the algebra $\Sym\Ee$ may depend, contrary to the
classical symmetries, on  arbitrary order derivatives. Next, in conformity with the above 
interpretation of the diffiety $(\Ee,\Cc^{\infty}_{\Ee})$, the informal interpretation of  
Definition (\ref{sym}) is : 
\begin{quote}
{\it Elements of the Lie algebra  $\Sym\Ee$ are \emph{vector \\ fields  on the ``space of all 
solutions of  $\Ee$\,".}}
\end{quote}
The importance of this interpretation is that it forces the question: 
\begin{quote}
{\it What are tensors, differential operators, PDEs, \\
etc. on the ``space of all solutions of  $\Ee$\,".}
\end{quote}
Later we shall give some examples and indications on how to define and use this kind of objects.
These objects form the thesaurus of {\it secondary calculus}, which is a natural language of the 
general theory of PDE's (see \cite{V2,KVe,KVin}). 

Finally, note that higher symmetries are not genuine vector fields as in the classical theory but
just some cosets of them modulo $\Cc_{\Ee}^{\infty}$. For this reason their action on functions
on $\Ee_{\infty}$ is not even defined. This at first glance discouraging fact leads to the 
bifurcation point: either to give up or to understand what are {\it functions on the ``space of 
solutions of  $\Ee$\,"}. Since, as we shall see, Definition (\ref{sym}), works well, the first
alternative should be discarded, while the second one will lead us to discover {\it differential
forms on the ``space of solutions of  $\Ee$\,"}.

\subsection{Computation of higher symmetries.} \label{comp-sym}
Though elements of $\vk$ are cosets of vector fields 
modulo $\Ee_{\infty}$ we can say that $\chi=[X]\in\vk$ is tangent to $\Ee_{\infty}$ if the vector
field $X$ is tangent to $\Ee_{\infty}$. Since $\Cc_{\infty}$ is tangent to $\Ee_{\infty}$, this 
definition is correct. If  $\Ee_{\infty}$ is locally given by equations (\ref{prolong}) and $\chi$
by the evolutionary derivation $\ev_{\varphi}$, then $\chi$ is tangent to $\Ee_{\infty}$ if and
only if $\ev_{\f}(D_{\sigma}(F_s))\left |\right._{\Ee_{\infty}}=0, \,\forall \sigma, s$. 
Since $\ev_{\f}$ and the $D_i$ commute these conditions are equivalent to 
$\ev_{\f}(F_s)\left |\right._{\Ee_{\infty}}=0, \,\forall s$, or, shortly, to  
$\ev_{\f}(F)\left |\right._{\Ee_{\infty}}=0$ with $F=(F_1,\dots,F_r)$. The bidifferential operator
$(\f,F)\mapsto \ev_{\f}(F)$ may be rewritten in the form $\ev_{\f}(F)=\ell_F(\f)$ with\newline 
\begin{equation}\label{lnz}
\ell_F=\left(\begin{array}{ccc}\sum_\sigma \frac{\partial F_1}{\partial u_\sigma^1}D_\sigma& \ldots & \sum_\sigma \frac{\partial F_1}{\partial u_\sigma^m}D_\sigma \\ \vdots &   & \vdots \\ \sum_\sigma \frac{\partial F_l}{\partial u_\sigma^1}D_\sigma & \ldots & \sum_\sigma \frac{\partial F_l}{\partial u_\sigma^m}D_\sigma\end{array}\right)
\end{equation}
and $\ell_F$ is called the {\it universal linearization operator}. Being tangent to  $\Ee_{\infty}$ 
the fields
$D_i$'s can be restricted to $\Ee_{\infty}$. It follows from (\ref{lnz}) that $\ell_F$  can also be 
restricted to $\Ee_{\infty}$. This restriction will be denoted by $\overline{\ell_F}$. So, by definition,
$\overline{\ell_F}(G\,|_{\Ee_{\infty}})=\ell_F(G)\,|_{\Ee_{\infty}}, \,\forall G$. In these terms
the condition of tangency of $\chi$ to $\Ee_{\infty}$ reads
\begin{equation}\label{eq.sym}
\boxed{\overline{\ell_F}(\bar{\f})=0, \;\;\bar{\f}=\f\,|_{\Ee_{\infty}}} \quad\Longrightarrow\quad
\boxed{\Sym\Ee=\ker\overline{\ell_F}}
\end{equation}
Hence the problem of the computation of the infinitesimal symmetries of a PDE $\Ee$ is reduced to
resolution of equation (\ref{eq.sym}). This equation is not a usual PDE, since it is imposed on
functions depending on unlimited number of variables. Nevertheless, it is not infrequent that
it can be exactly solved. For instance, this method allows not only to easily rediscover ``classical" hierarchies associated with well known integrable systems but also to find various 
new ones (see \cite{KVin,KK}).

The interpretation of higher symmetries as vector fields on the ``space of solutions of  
$\Ee$\," leads to the question : What are the trajectories of this field? The equation of trajectories 
of $\chi$ is very natural:
\begin{equation}\label{eq.traec}
u_t=\varphi(x,u,\ldots, u_\sigma^i,\ldots) \quad\mbox{with} \quad \f=(\f_1,\dots,\f_m).
\end{equation}
Equation (\ref{eq.traec}) is the exact analogue of the classical equations 
\begin{equation}\label{cl.traec}
x_i=a_i(x), \;\;i=1,\dots,m, \;\;x=(x_1,\dots,x_m),
\end{equation}
which describe trajectories of the vector field $X=\sum_{i=1}^{m}a_i\p/\p x_i$. An essential 
difference between  equations (\ref{eq.traec}) and  (\ref{cl.traec}) is that the initial data uniquely 
determine solutions of  (\ref{cl.traec}), while it is not longer so for  (\ref{eq.traec}). Indeed, 
the uniqueness for the partial evolution equation is guaranteed by some additional to the initial  
conditions, for instance, the boundary ones. For this reason a ``vector field" $\chi\in\vk$  does 
not generate a flow on the ``space of solutions of  $\Ee$\,". 

A very important consequence of this fact is that in this new context the classical relation
between Lie algebras and Lie groups breaks down. Consequently, the absolute priority  
should be given to infinitesimal symmetries, not to the finite ones. 

One of the most popular applications of symmetry theory  takes an especially simple
form if expressed in terms of generating functions. Namely, imagine for a while that the flow
generated by $\chi\in\vk$ exists. Then, according to (\ref{eq.traec}), ``stable points" of this
flow are solutions of the equation $\f=0$. In other words, these ``stable points" are solutions 
of the last equation. If $\f^1,\dots,\f^l$ are generating functions of some symmetries of $\Ee$,
then solutions of the system
\begin{equation}\label{inv-sol}
\left\{\begin{array}{c}F=0 \\\varphi^1=0 \\\vdots \\\varphi^l=0\end{array}\right.
\end{equation}
represent those solutions of $\Ee$ that are stable in the above sense with respect to ``flows"
generated by $\f^1,\dots,\f^l$. System (\ref{inv-sol}) is well overdetermined and by this reason
can be exactly solved in many cases. For instance, famous multi-soliton solutions of the KdV 
equation are solutions of this kind.

\subsection{What are partial differential equations?} The fact that we have built a self-consistent
and well working theory of symmetries for PDEs based on diffieties gives a considerable
reason to recognize diffieties as objects of category of PDEs. Another argument supporting
this idea is as follows.

Take any PDE, say,\newline
\begin{equation}\label{eq1}
u_{xx}u_{tt}^2+u_{tx}^2+(u_x^2-u_t)u=0.
\end{equation}
This is a hypersurface  $\Ee\subset J^2(E,2), \,\dim E=3$. The equivalent system of first
order PDEs is
\begin{equation}\label{eq2}
\left\{\begin{array}{l}u_x=v \\u_t=w \\v_xw_t^2+v_tw_x+(v^2-w)u=0 .\end{array}\right.\ 
\end{equation}
This is a submanifold  $\Ee'\subset J^1(E',2), \,\dim E'=5,$ of codimension 3. $\Ee$ and 
$\Ee'$ live in different jet spaces and have different dimensions. For this reason their 
classical symmetries cannot even be compared. On the other hand, associated with
$\Ee$ and $\Ee'$ diffieties are naturally identified and hence have the same (higher)
symmetries. So, this fact may be interpreted by saying that (\ref{eq1}) and (\ref{eq2})
are different descriptions of the same object, namely, of the associated diffiety.

Another example illustrating priority of diffieties is the factorization problem. Namely, if
$G$ is a Lie algebra of classical symmetries of an equation $\Ee$, then the question is:
Can  $\Ee$ be factorized by the action of $G$ and what is the resulting ``quotient 
equation"\,?  In terms of diffieties the answer is almost obvious : this is the equation
$\Ee'$ such that $\Ee_{\infty}\backslash G=\Ee'_{\infty}$. On the contrary, it is not very 
clear how to answer this question in terms of the usual approach. 
\begin{ex}\label{fact-L}
Let $G$ be the group of translations of the Euclidean plane. Obviously, these translations 
are symmetries of the Laplace equation $u_{xx}+u_{yy}=0$. Then the corresponding quotient
equation is again the Laplace equation.
\end{ex}
There are many other examples manifesting that
\begin{quote}
{\it A PDE as a mathematical object is a \emph{diffiety}, while what is usually called a PDE is
just one of many possible ``identity cards" of it.}
\end{quote}
It should be stressed that the diffiety associated with a system of PDEs (in the usual sense 
of this word) is the exact analogue of the algebraic variety associated to a system of algebraic 
equations. Indeed, if a system of algebraic equations is $f_1=0,\dots,f_r=0$, then the ideal
defining the corresponding variety is algebraically generated by polynomials 
$f_i$ In the case of a PDE $\Ee=\{F_i=0\}$ the ideal defining $\Ee_{\infty}$ is algebraically generated not only by functions $F_i$ but also by all their differential consequences 
$D_{\sigma}(F_i)$ (see (\ref{prolong-coord})). Viewed from this side algebraic geometry is seen
as the zero-dimensional case of the general theory of PDEs.

\section{On the internal structure of diffieties.}\label{inter-dif} 
On the surface, a diffiety $\Oo=(\Ee,\Cc^{\infty}_{\Ee})$ looks  as a simple enough object 
like a foliation. All foliations of given finite dimension and
codimension are locally equivalent. On the contrary, the situation drastically changes when the
codimension becomes infinite. So, the problem of how to extract all the information on the equation
$\Ee$, which is encoded in the ``poor" Frobenius distribution $\Cc^{\infty}_{\Ee}$, naturally arises
and becomes central. To gain a first insight into the problem we consider as a simple model a 
Frobenius distribution $\Dd$, or, equivalently, a foliation, on a finite-dimensional manifold $M$. 

\subsection{The normal complex of a Frobenius distribution.}\label{Norm-Frob}
Let $\Dd$ be an $r$--dimensional Frobenius distribution on a manifold $M$. 
The quotient $\ci(M)$--module 
$\Nn=D(M)/\Dd$ is canonically isomorphic to $\G(\nu)$ where 
$\nu$ is the normal to $\Dd$ bundle, i.e., the bundle whose fiber over $x\in M$ is 
$T_xM/\Dd_x$. Put $\hat{Y}= [Y\!\mod \Dd ]\in \Nn$ for  $Y\in D(M)$ and
$$
\nabla_X(\hat{Y})=\widehat{[X,Y]} \quad\mbox{for} \quad X\in\Dd.
$$  
It is easy to see that $\nabla_{fX}=f\nabla_X, \,\nabla_X(f\hat{Y})=X(f)\hat{Y}+f\nabla_X(\hat{Y})$ 
if $f\in \ci(M)$ and $[\nabla_X,\nabla_{X'}]=\nabla_{[X,X']}$.  These formulas tell that the correspondence $\nabla:X\mapsto \nabla_X$ is a flat {\it $\Dd$--connection}. This means that 
this construction can be restricted to a leaf $\Ll$ of the foliation associated with $\Dd$ and this 
restriction is a flat connection $\nabla^{\Ll}$ in the normal to $\Dd$ bundle $\nu$ restricted to 
$\Ll$. Recall that with a flat connection s associated a de Rham-like complex i (see \cite{DPV}), 
which for $\nabla^\Ll$ is
\begin{equation}\label{dRcon}
0\lrw\Nn_{\Ll}\stackrel{\nabla^{\Ll}}{\lrw}\Lambda^1(\Ll)\otimes_{\ci(\Ll)}\Nn_{\Ll}
\stackrel{\nabla^{\Ll}}{\lrw}\dots
\stackrel{\nabla^{\Ll}}{\lrw}\Lambda^r(\Ll)\otimes_{\ci(\Ll)}\Nn_{\Ll}\lrw0
\end{equation}
where the covariant differential is abusively denoted also by $\nabla^{\Ll}$ and 
$\Nn_{\Ll}=\Gamma(\nu|_{\Ll})$.
This complex is, in fact, the restriction of the complex to $\Ll$.
\begin{equation}\label{dRD}
0\lrw\Nn\stackrel{\nabla}{\lrw}\Lambda^1_{\Dd}\otimes_{\ci(M)}\Nn
\stackrel{\nabla}{\lrw}\dots\stackrel{\nabla}{\lrw}\Lambda^r_{\Dd}\otimes_{\ci(M)}\Nn\lrw0
\end{equation}
where $\Lambda^i_{\Dd}=\Lambda^i( M)/\Dd\Lambda^i(M)$ with 
$$
\Dd\Lambda^i(M)=\{\omega\in\Lambda^i(M)\,|\,\omega(X_1,\dots,X_i)=0,  
\;\forall X_1,\dots,X_i\in\Dd \}.
$$
The terms of the complex (\ref{dRD}) are $\Nn$--valued differential forms on $\Dd$, i.e.,
$\rho(X_1,\dots,X_s)\\\in\Nn$ if $X_1,\dots,X_s\in\Dd$. The {\it covariant differential} $\nabla$
is defined as
$$
\begin{array}{r}
\nabla(\rho)(X_1,\dots,X_{s+1})=
\sum_{i=1}^{s+1}(-1)^{i-1}\nabla_{X_i}(\rho(X_1,\dots,\widehat{X_i},\dots,X_{s+1}))+\bigskip\\
\sum_{i<j}(-1)^{i+j}\rho([X_i,X_j],X_1\dots,X_i,\dots,X_j,\dots,X_{s+1}).
\end{array}
$$ 
Pictorially, this situation may be seen as a ``foliation" of the complex (\ref{dRD}) by complexes
(\ref{dRcon}). The i-th cohomology of complexes (\ref{dRcon}) and (\ref{dRD}) will be denoted by $H^i(\nabla^{\Ll})$  and $H^i(\nabla)$, respectively.  We also have a natural restriction map $H^i(\nabla)\rw H^i(\nabla^{\Ll})$ in cohomology.

Formally, the above construction remains valid for any Frobenius distribution and hence can be applied to diffieties. In order to duly specify complex (\ref{dRD}) to this particular case we need
a new construction from differential calculus over commutative algebras.

\subsection{Modules of jets.} Let $A$ be an unitary algebra and let $P,Q$ be $A$--modules.
Denote by $\Diff_k(P,Q)$ the totality of DO's of order $\leq k$ considered as a left
$A$--module, i.e., $(a,\square)\mapsto a\square, \,a\in A, \,\square\in\Diff_k(P,Q)$. 
Consider a subcategory  $\Kk$  of the category of $A$-modules such that 
$\Diff_k(P,Q)\in \mathrm{Ob \,\Kk}$ if $P,Q\in \mathrm{Ob \,\Kk}$. For a fixed $P$ we have
the functor $Q\mapsto \Diff_k(P,Q)$. We say that a pair composed of an $A$--module
$\Jj^k_{\Kk}(P)$ and a $k$-th order DO $j_k=j_k^{P,\Kk}:P\rw \Jj^k_{\Kk}$ represents this 
functor in the category $\Kk$ if the map 
$\Hom_A(\Jj^k_{\Kk}(P),Q)\ni h\mapsto h\circ j_k\in \Diff_k(P,Q)$
is an isomorphism of $A$--modules. Under some weak condition on $\Kk$, 
which we skip, the {\it representing object} $(\Jj^k_{\Kk}(P),j_k)$ exists and is unique 
up to isomorphism. $\Jj^k_{\Kk}(P)$ is called the module of {\it k-th order jets of $P$} 
(in $\Kk$). Thus for a DO $\square\in\Diff_k(P,Q)$ there is a unique $A$--module 
homomorphism $h_{\square}:\Jj^k_{\Kk}(P)\rw Q$ such that 
$\square=h_{\square}\circ j_k$. 

As an 
$A$--module $\Jj^k_{\Kk}(P)$ is generated by elements $j_k(p), \,p\in P$. A natural transformation 
of functors $\Diff_l(P,\cdot)\mapsto\Diff_k(P,\cdot), \,l\leq k,$ induces a homomorphism
$\pi_{k,l}=\pi_{k,l}^P:\Jj^k_{\Kk}(P)\rw\Jj^l_{\Kk}(P)$ of $A$--modules such that 
$j_l=\pi_{k,l}\circ j_k$. This allows to define the inverse limit of pairs
$(\Jj^k_{\Kk}(P),j_k)$ called the module of {\it infinite jets} of $P$ and denoted by 
$(\Jj^{\infty}_{\Kk}(P), \,j_{\infty}=j_{\infty}^{P,\Kk})$. Natural projections 
$\pi_{\infty,k}:\Jj^{\infty}_{\Kk}(P)\rw\Jj^k_{\Kk}(P)$ come from the definition. These maps supply 
$\Jj^{\infty}_{\Kk}(P)$ with a decreasing filtration
\begin{equation}\label{j-filt}
\Jj^{\infty}_{\Kk}(P)\supset\ker(\pi_{\infty,0})\supset\ker(\pi_{\infty,1})\supset\dots
\supset\ker(\pi_{\infty,k})\supset\cdots
\end{equation}
Finally, we stress that 
$\Jj^k_{\Kk}(P)$ and all related constructions essentially depend on $\Kk$.

Any operator $\square\in \Diff_r(P,Q)$ induces a homomorphism 
$$
h_{\square}^r : \Jj^{k+r}_{\Kk}(P)\rw\Jj^r_{\Kk}(Q), \quad r\geq 0.
$$ 
Namely, the composition 
$P\stackrel{\square}{\lrw}Q\stackrel{j_r}{\lrw}\Jj^r_{\Kk}(Q))$ is a DO of order $\leq k+r$.  
So, it can be presented in the form $h_{j_r\circ\square}\circ j_{k+r}$, and we put
\begin{equation}\label{h-square}
h_{\square}^r\df h_{j_r\circ\square} : \Jj^{k+r}_{\Kk}(P)\rightarrow \Jj^r_{\Kk}(Q).
\end{equation}
The inverse limit of homomorphisms $h_{\square}^r$'s defines a homomorphism
of filtered modules
$$
h_{\square}^{\infty} : \Jj^{\infty}_{\Kk}(P)\rightarrow \Jj^{\infty}_{\Kk}(Q)
$$
which shifts filtration (\ref{j-filt}) by $-k$.

If $Q=\Jj^{k}_{\Kk}(P)$ and $\square=j_k$, then the above construction gives natural
inclusions
$$
\iota_{k,r}\df h_{j_r\circ j_k}^r : \Jj^{k+r}_{\Kk}(P)\hookrightarrow \Jj^r_{\Kk}(\Jj^{k}_{\Kk}(P)).
$$
Their inverse limit of these inclusions is
$$
\iota_{\infty} : \Jj^{\infty} _{\Kk}(P)\hookrightarrow \Jj^{\infty} _{\Kk}(\Jj^{\infty}_{\Kk}(P)).
$$

Now we shall describe constructively the above {\it conceptually} defined modules
of jets for {\it geometrical modules} over the algebra $A=\ci(M)$. Recall that an $A$--module
$P$ is geometrical if all its elements $p$ such that $p\in\mu_z\!\cdot \!P, \,\forall z\in M,$ are
equal to zero (see \cite{Nestr}). Here $\mu_z=\{f\in\ci(M)\,|\,f(z)=0\}$. The category of
geometrical $A$--modules will be denoted by $\mathcal{G}$ and we shall write simply 
$\Jj^k(P)$ for $\Jj^k_{\mathcal{G}}(P)$. \newline
 
Put $\Jj^k=\Jj^k(A)$ and note that $\Jj^k$ is a unitary algebra with the product
$(f_1j_k(g_1))\cdot (f_2j_k(g_2))=f_1f_2j_k(g_1g_2), \,f_i, g_i\in A$. In particular, $\Jj^k$
is a bimodule. Namely, left (standard) and right multiplications by $f\in A$ are defined as
$(f,\theta)\mapsto f\theta$ and $(f,\theta)\mapsto \theta j_k(f)$, respectively. $\Jj^k$
supplied with the right $A$-module structure will be denoted by $\Jj^k_{>}$.

We have (see \cite{Nestr}).
\begin{proc}\label{jets}
\begin{enumerate}
\item Let $\a_k$ be the vector bundle whose fiber over $z\in M$ is $J_z(M)$ (see 
section \ref{JS}). Then $\Jj^k=\Gamma(\a_k)$.
\item $\Jj^k(P)=\Jj^k_{>}\otimes_AP$.
\end{enumerate}
\end{proc}

\subsection{Jet-Spencer complexes.}
The $k$--th jet-Spencer complex of $P$ denoted by $\mathcal{S}_k(P)$ is defined as
\begin{equation}\label{j-Sp}
\begin{array}{rl}
0\rw \Jj^k(P)\stackrel{S_k}{\lrw}\Jj^{k-1}(P)\otimes_A\Lambda^1(M)
\stackrel{S_k}{\lrw}\dots  
\stackrel{S_k}{\lrw}\Jj^{k-n}(P)\otimes_A\Lambda^n(M)\rw 0
\end{array}
\end{equation}
with $S_k(j_{k-s}(p)\otimes\omega)=j_{k-s-1}(p)\otimes d\omega, \,\omega\in\Lambda^s(M)$.
Here $n=\dim M$ and we assume that $\Jj^s(P)=0$ if $s<0$.

Differentials of Spencer complexes are 1-st order DOs. For $k\geq l$, the homomorphisms
$$
\pi_{k-s,l-s}\otimes\id|_{\Lambda^s(M)}:\Jj^{k-s}(P)\otimes_A\Lambda^s(M)\lrw
\Jj^{l-s}(P)\otimes_A\Lambda^s(M), \,s=0,1,\dots,n,
$$
define a cochain map $\sigma_{k,l}:\mathcal{S}_k(P)\lrw\mathcal{S}_l(P)$. In particular,
we have the following sequence of Spencer complexes
\begin{equation}\label{j-Sp-tower}
0\lw\Ss_0(P)\stackrel{\sigma_{1,0}}{\llw}\Ss_1(P)\stackrel{\sigma_{2,1}}{\llw}\dots
\stackrel{\sigma_{k,k-1}}{\llw}\Ss_k(P)\stackrel{\sigma_{k+1,k}}{\llw}\dots
\end{equation}
The infinite jet-Spencer complex $\Ss_{\infty}(P)$ is defined as the inverse limit of 
(\ref{j-Sp-tower}) together with natural cochain maps $\sigma_{\infty,k}$. As in the case of jets
the complex $\mathcal{S}_{\infty}(P)$ is filtered by subcomplexes $\ker(\sigma_{\infty,k})$.

An operator $\square\in \Diff_r(P,Q)$ induces a cochain map of Spencer complexes
\begin{equation}\label{map-Sp}
\s_{\square}^{k} : \mathcal{S}_k(P)\lrw\mathcal{S}_{k-r}(Q) ,
\end{equation}
which acts on the $s$--th term of $\Ss_k(P)$ as
\begin{equation}\label{map-Sp-term}
h_{\square}^{k-s}\otimes\id|_{\Lambda^s(M)}:\Jj^{k-s}(P)\otimes_A\Lambda^s(M)\lrw
\Jj^{k-s-r}(Q)\otimes_A\Lambda^s(M)
\end{equation}
(see (\ref{h-square})).

Conceptually, the $k$--th jet-Spencer complex is an {\it acyclic resolvent} for the 
{\it universal} $k$--order differential operator $j_k$. Namely, we have
\begin{proc}\label{Sp-coho}
If $P=\G(\xi)$ with $\xi$ being a vector bundle over $M$, then 
\begin{enumerate}
\item $\Ss_k(P)$ is acyclic in positive dimensions $\Leftrightarrow \, H^{i}(S_k)=0$
if $i>0$.
\item $H^{0}(S_k)=P$ and $0$--cocycles  are $j_k(p)\in \Jj^k(P), \,p\in P$.
\end{enumerate}
\end{proc} 

Jet-Spencer complexes are {\it natural}, since they can be defined over arbitrary unitary 
(graded) algebras (see \cite{V4}). In other words, they are \emph{compatible} with 
homomorphisms of these algebras. In particular, they restrict to submanifolds. For our 
purposes, we need to describe this procedure.

Let $N\subset M$ be a submanifold. In the notation of Proposition \ref{Sp-coho} we
put $P_N=\G(\xi|_N)$ and $j_k^N:P_N\rw \Jj^k(P)$ to distinguish this jet-operator
on $N$ from $j_k:P\rw \Jj^k(P)$. Since $\ci(N)$--modules can also be considered as 
$\ci(M)$--modules, the composition $\square$
$$
P\stackrel{restriction}{\lrw} P_N\stackrel{j_k^N}{\lrw} \Jj^k(P_N)
$$
is a $k$--th order DO over $\ci(M)$. The homomorphism of $\ci(M)$--modules 
$h^{\square}:\Jj^k(P)\rw \Jj^k(P_N)$ associated with $\square$ is, by definition, the 
{\it restriction operator}. Now, by tensoring this restriction operator with the well-known
restriction operator for differential forms we get the restriction operator for terms of
$\Ss_k(P)$. Finally, by passing to the inverse limit we get the restriction operator for
$\Ss_{\infty}(P)$.

\subsection{Foliation of Spencer complexes by a Frobenius distribution.}\label{hor}
Take the notation of Subsection \ref{Norm-Frob} and denote by $\Dd\Lambda^i(M)$
(resp., $\Dd\Jj^k(P))$ the totality of all differential forms (resp., jets) whose restrictions
to all leaves of $\Dd$ are trivial. Similarly,  $\Dd \Ss_k(P)$ stands for the maximal 
subcomplex of $\Ss_k(P)$ such that its restrictions to leaves of $\Dd$ are trivial.
{\it Horizontal} (with respect to $\Dd$) differential forms and jets are elements of
the quotient modules
$$ 
\bar{\Lambda}^i_{\Dd}(M)\df\Lambda^i(M)/\Dd\Lambda^i(M), \quad
\bar{\Jj}^k(P)\df\Jj^k(P)/\Dd\Jj^k(P) .
$$
Similarly, the {\it horizontal jet-Spencer complex} is
$$
\bar{\Ss}_k(P)\df S_k(P)/\Dd \Ss_k(P) .
$$
Restrictions of all the above horizontal objects to leaves of $\Dd$ are naturally defined.
Conversely, a horizontal differential form (resp., jet, Spencer complex) may be viewed
as a family of differential forms (resp., jet, Spencer complex) defined on single leaves.
In other words, if $\Ll$ runs the all leaves of $\Dd$, then the $\Lambda^i(\Ll)$ foliate 
$\bar{\Lambda}^i_{\Dd}(M)$, and similarly for jets and Spencer complexes.

Accordingly, the exterior differential $d$ as well as the Spencer differential $S_k$ 
factorize to $\bar{\Lambda}^*(M)$ and $\bar{\Ss}_k(P)$, since   $\Dd\Lambda^i(M)$
and  $\Dd \Ss_k(P)$ are stable with respect to $d$ and $S_k$, respectively. These 
quotient differentials will be denoted by $\bar{d}$ and $\bar{S}_k$, respectively.
In this way we get the {\it horizontal de Rham} and {\it Spencer complexes} and 
hence {\it horizontal de Rham} and {\it Spencer cohomology}.
\begin{rmk}\label{hor-dir}
For simplicity, in the above definition of horizontal objects we have used leaves of $\Dd$.
It fact, with a longer formal procedure this can be done explicitly in terms of the
distribution $\Dd$.
\end{rmk}

\subsection{The normal complex of a diffiety and more symmetries.} 
First, we shall describe the normal complex for 
$J^{\infty}(E,n)$. Denote the normal bundle to $\Cc^k$ 
by $\nu_k, \,1\leq k\leq \infty$ and put $\vk_k\df\Gamma(\nu_k)=D(J^k(E,n))/\Cc^k$, 
$\vk_{l,k}\df\Gamma(\pi_{k,l}^*(\nu_l))=\pi_{k,l}^*({\vk_l}), \,l\leq k$.
\begin{proc}\label{J-coho}
\begin{enumerate}
\item $\vk=\vk_{1,\infty}, \;\vk_{k,\infty}=\bar{\Jj}^k(\vk)$ and
$\vk_{\infty}=\bar{\Jj}^{\infty}(\vk)$;
\item the normal to $\Cc^{\infty}$ complex is isomorphic to $\bar{\Ss}_{\infty}(\vk)$;
\item $H^0(\bar{S}_{\infty})=\vk, \;H^i(\bar{S}_{\infty})=0, \;i\neq  0$.
\end{enumerate}
\end{proc}
Assertion (3) in this proposition is a pro-finite consequence of Proposition \ref{Sp-coho}.
The following important interpretation is a consequence of this and the second assertions
of the Proposition \ref{J-coho}: 
\begin{quote} 
{\it $\vk$ is the zero-th cohomology space of the  complex \\ normal  to the infinite 
contact structure $\Cc^{\infty}$.}
\end{quote}

Now we can describe the normal complex to $\Cc^{\infty}_{\Ee}$  by restricting, in a sense, Proposition \ref{J-coho} to $\Ee_{\infty}$.  First, to this end, we need a conceptually
satisfactory definition of the universal linearization operator (\ref{lnz}), which was defined 
coordinate-wisely in Subsection \ref{comp-sym}. 

The equation $\Ee\subset J^k(E,n)$ may be presented in a coordinate-free form as 
$\Phi=0, \,\Phi\in \G(\xi)$ with $\xi$ being a suitable vector 
bundle over $J^k(E,n)$. If $P=\G(\pi_{\infty,k}^*(\xi))$ and $\Phi_{\infty}=\pi_{\infty,k}^*(\Phi)$, 
then $\Ee_{\infty}=\{\bar{j}_{\infty}(\Phi_{\infty})=0\}$. Additionally, assume that $P$ is supplied 
with a connection $\nabla$. If $Y\in\Cc^{\infty}$, then, as it is easy to see, 
$\nabla_Y(\Phi)|\Ee_{\infty}=0$. For this reason the following definition is correct.
\begin{equation}
\ell_{\Ee}(\chi)\df \nabla_{X}(\Phi)|\Ee_{\infty} \quad \mbox{with} \quad \chi=[X\!\!\!\mod \Cc^{\infty}_{\Ee}].
\end{equation}
The operator $\ell_{\Ee}:\vk\rw P$ does not depend on the choices of $\Phi$ and $\nabla$. It defines
the cochain map of jet-Spencer complexes
\begin{equation}
\sigma_{\ell_{\Ee}}^{\infty}:\Ss_{\infty}(\vk)\lrw \Ss_{\infty}(P)
\end{equation}
(see (\ref{map-Sp})).
\begin{proc}\label{sym-grad}
Let $\Nn_{\Ee}$ be the normal to the distribution $\Cc^{\infty}_{\Ee}$ complex on 
$\Ee_{\infty}$. Then
\begin{enumerate}
\item $\Nn_{\Ee}$ is isomorphic to the complex $\ker \sigma_{\ell_{\Ee}}^{\infty}$. 
\item The cohomology $H^i(\Nn_{\Ee})$ of the complex  $\Nn_{\Ee}$ is trivial if $i>n$.
\item $\Sym\,\Ee=H^0(\Nn_{\Ee})=\ker\,\ell_{\Ee}$.
\end{enumerate}
\end{proc}
Assertions (2) and (3) of this proposition are consequences of the first one, which allows
to compute the cohomology of  $\Nn_{\Ee}$. Moreover,
assertion (3) is one of many other arguments that motivate the following definition.
\begin{defi}\label{def-sym}
The \emph{Lie algebra of (higher) infinitesimal symmetries of a PDE $\Ee$} is the 
cohomology of the normal complex  $\Nn_{\Ee}$.
\end{defi}
Accordingly, denote by $\Sym_i\Ee$ the $i$--th cohomology space of  $\Nn_{\Ee}$. So, the
whole Lie algebra of infinitesimal symmetries of $\Ee$ is {\it graded}\,:
$$
\boxed{\Sym_{*}\Ee=\sum_{i=0}^{n}\Sym_i\Ee, \quad \Sym_i\Ee=H^i(\Nn_{\Ee})}
$$
In particular, $\Sym\Ee =\Sym_0\Ee$ (see Subsection \ref{High-sym}). 
\begin{rmk}
The description of the Lie product in $\Sym_{*}\Ee$ is not immediate and requires some new 
instruments of differential calculus over commutative algebras. For this reason we shall
skip it.
\end{rmk}
The following proposition illustrates what the algebra  $\Sym_{*}\Ee$ looks like.
\begin{proc}
If $\Ee$ is not an overdetermined system of PDEs, then 
$$
\Sym_{0}\Ee=\ker\ell_{\Ee}, \quad \Sym_{1}\Ee=\mathrm{coker} \ell_{\Ee}
\quad\mbox{and} \quad\Sym_{i}\Ee=0\quad \mbox{if}\quad i\neq 0,1. 
$$
\end{proc}
In this connection we note that a great majority of the PDEs of current interest in geometry,
physics and mechanics are not overdetermined. As an exception we mention the system of
Yang-Mills equations, which is sightly overdetermined, and  for these equations 
$\Sym_{2}\Ee\neq 0$.
 
 To conclude this section we would like to emphasize the role of the structure of differential 
 calculus over commutative algebras in the above discussion. While we have used the 
 ``experimental data" coming from the theory of integrable systems  to discover ``by hands" 
 the conceptually simplest part of infinitesimal symmetries of PDEs, i.e., the Lie algebra 
 $\Sym\Ee$, a familiarity with the structures of differential calculus over 
 commutative algebras is indispensable to discover that it is just the zeroth component of
 the full symmetry algebra $\Sym_{*}\Ee$, which in its turn is the cohomology of a certain 
 complex. 
 
 \section{Nonlocal symmetries and once again : what are PDEs ?}
 In the previous section we have constructed a self-consistent symmetry theory, which,
 from one side, resolves shortcomings of the classical theory discussed in Sections 2-5
 and, from another side, incorporates ``experimental data" that emerged in the theory 
 of integrable systems. However, one important element of this theory was taken into
 account. Namely, we have in mind {\it nonlocal symmetries}. Roughly
 speaking, these are symmetries whose generating function depends on variables of the
 form $D_i^{-1}(u)$. Fortunately, these unusual symmetries can be tamed by introducing 
 only one new notion we are going to describe.

\subsection{Coverings of a diffiety.} 
Schematically, a diffiety $\Oo$ is a {\it pro-finite manifold}  $\M$ supplied with a 
finite-dimensional pro-finite Frobenius distribution $\Dd=\Dd_{\Oo}$\,:\,$\Oo=(\M,\Dd)$. 
We omit technical details that these data must satisfy.

Recall that a pro-finite manifold is the inverse limit of a sequence of smooth maps
\begin{equation}\label{pro-mnf}
M_0\stackrel{\mu_1}{\llw}M_1\stackrel{\mu_2}{\llw}\dots \stackrel{\mu_k}{\llw}
M_k\stackrel{\mu_{k+1}}{\llw}\dots \quad \Leftarrow\M
\end{equation}
A pro-finite distribution on $\M$ is the inverse limit via $\mu_i$'s of distributions $\Dd_i$'s 
on $M_i$'s. The associated sequence of homomorphisms of smooth function algebras
\begin{equation}\label{pro-alg}
\ci(M_0)\stackrel{\mu_1^*}{\lrw}\ci(M_1)\stackrel{\mu_2^*}{\lrw}\dots \stackrel{\mu_k^*}{\lrw}
\ci(M_k)\stackrel{\mu_{k+1}^*}{\lrw}\dots \quad \Rightarrow \Ff_{\M}
\end{equation}
with $\Ff_{\M}$ being the direct limit of homomorphisms $\mu_k^*$ is filtered by subalgebras
$\Ff_{\M}^k\df\mu_{\infty,k}^*\ci(M_k)$ where $\mu_{\infty,k}:\M\rw M_k$ is a natural projection. Differential calculus on $\M$ is interpreted as the calculus over the filtered algebra $\Ff_{\M}$
(see Subsection \ref{dif-cal}). The dimension of $\Dd$ is interpreted as the ``number of 
independent variables".

A morphism $F:\Oo\rw\Oo^{\prime}$ of a diffiety $\Oo=(\M,\Dd)$ to a diffiety 
$\Oo^{\prime}=(\M^{\prime},\Dd^{\prime})$ is, abusing the notation, a map $F:\M\rw\M^{\prime}$
such that $F^*(\Ff_{\M^{\prime}})\subset \Ff_{\M}, \,F^*$ is compatible with filtrations and 
$d_{\theta}F(\Dd_{\theta})\subset \Dd_{F(\theta)}^{\prime}, \,\forall \theta\in\M$. 
\begin{defi}\label{cover}
A surjective morphism $F:\Oo\rw\Oo^{\prime}$ of diffieties is called a \emph{covering} if 
$\dim \Dd=\dim \Dd^{\prime}$ and $d_{\theta}F$ isomorphically sends $\Dd_{\theta}$ to
$\Dd_{F(\theta)}^{\prime}, \,\forall \theta\in\M$.
\end{defi}
This terminology emphasizes the analogy with the standard notion of a covering in the category
of manifolds. Namely, fibers of a covering are zero-dimensional diffieties in the sense that the
their structure distributions $\Dd$'s are zero-dimensional. If these fibers are finite-dimensional
in the usual sense, then the covering is called {\it finite-dimensional}. 

A covering $F:\Ee_{\infty}\rw\Ee_{\infty}^{\prime}$ may be interpreted as a (nonlinear) DO, 
which sends solutions of $\Ee$ to solutions of $\Ee^{\prime}$. More exactly, it associates 
with a solution of $\Ee^{\prime}$ a families of solutions of $\Ee$. For instance, the famous
Cole-Hopf substitution $v=2u_x/u$ that sends solutions of the heat equation 
$\Ee=\{u_t=u_{xx}\}$ to solutions of the Burgers equation $\Ee^{\prime}=\{v_t=v_{xx}+v v_x\}$ 
comes from a 1-dimensional covering of $\Ee^{\prime}$. Equivalently, the passage from
a PDE to a covering equation is the {\it inversion of a $($nonlinear$)$ DO} on solutions of this PDE.
For instance, by inverting the operator $v\mapsto 2v_x/v$ on solutions of the Burgers equation
one gets the heat equation.

\subsection{Where coverings  appear.}\label{ex-cover} 
The notion of a covering of a diffiety was introduced by the author (see \cite{Vcat}) as a
common basis for various constructions that appeared in PDE's. Below we list and briefly 
discuss some of them. 

1) In the language of diffieties the passage from Lagrange's description of a continuum media 
to that of Euler is interpreted as a covering. This interpretation allows to apply instruments of
secondary calculus to this situation and, as a result, to derive from this fact some  important consequences for mechanics of continua. \newline 

2) {\sl \underline{Factorization of PDE's.}} If $G$ is a symmetry group of a diffiety $\Oo$,  then 
under some natural conditions the quotient diffiety $\Oo\backslash G$ is well-defined and
$\Oo\rw \Oo\backslash G$ is a covering. In particular, if $\Oo=\Ee_{\infty}$, then 
$\Oo\backslash G=\Ee_{\infty}^{\prime}$. In such a case $\Ee^{\prime}$ is the
{\it quotient equation} of $\Ee$ by $G$. A remarkable fact is that the group $G$ 
in this construction may be an ``infinite-dimensional" Lie group like the group 
$\mathrm{Diffeo}(M)$ of diffeomorphisms of a manifold $M$, or the group of 
contact transformations, etc. \newline

3) {\sl \underline{Differential invariants and characteristic classes}} Let $\pi:E\rw M$ be a fiber 
bundle of geometrical structures of a type $\gS$ on $M$ 
(see \cite{ALV}). Then $\mathrm{Char}\,\gS=J^{\infty}(\pi)\backslash \mathrm{Diffeo}(M)$ 
is the {\it characteristic diffiety} for $\gS$-- structures. This diffiety is with singularities, 
which are in turn diffieties with a smaller numbers of independent variables. Functions
on $\mathrm{Char}\,\gS$ are scalar differential invariants of $\gS$-- structures, horizontal
de Rham cohomology is composed of their characteristic classes, etc.

Similarly one can define differential invariants and characteristic classes for solutions 
of {\it natural} PDEs, i.e., those that are invariant with respect to the group 
$\mathrm{Diffeo}(M)$ or some more specific subgroups of this group. For instance, 
Einstein equations and many other equations of mathematical physics are natural. 
Gel'fand-Fuks characteristic classes are quantities of this kind. The reader will find 
more details and examples in \cite{Vscal,MVY,Tsu}.\newline
 
4) {\sl \underline{B\"acklund transformations.}}  The notion of covering allows to rigorously
define B\"acklund transformations. Namely, the diagram
\begin{equation*}\label{Backl}
\xymatrix{&\Ee_{\infty}\ar[dl]_{F'}\ar[dr]^{F''}\\
\Ee_{\infty}^{\prime}&\qquad &\Ee_{\infty}^{''}}
\end{equation*}
where $F'$ and $F''$ are coverings presents the B\"acklund transformation $F''\circ(F')^{-1}$
from $\Ee'$ to $\Ee''$ and its inverse $F'\circ(F'')^{-1}$. The importance of this definition lies 
in the fact that it suggests an efficient and regular method for finding  B\"acklund transformations 
for a given PDE (see \cite{KVN,KVin,IgoB}). Previously this was a kind of handcraft art.
Moreover, it turned out to be possible to prove for the first time nonexistence of 
B\"acklund transformations connecting two given equations (see \cite{Igo}). This seems to 
be an impossible task by using only the standard  techniques of the theory of integrable systems.
\newline

5) {\sl \underline{Poisson structures and the Darboux lemma in field theory}.} 
The efficiency and elegance of the Hamiltonian approach to the mechanics of systems with 
a finite number of degrees of freedom motivates to look for its extension to the mechanics 
of continua and field theory.  Obviously, this presupposes a due
formalization of the idea of a Poisson structure in the corresponding infinite-dimensional 
context.  Over the past 70-80 years various concrete constructions of the Poisson bracket 
in field theory were proposed, mainly,  by physicists. But the first attempts to build a 
systematic general theory can be traced back only to late 1970's. Here we mention B.~A.~Kuperschmidt's paper \cite{Kuper} where he constructs an analogue of the Poisson 
structure on the cotangent bundle on infinite jets, and the paper by I.~M.~Gel'fand and 
I.~Dorfman  \cite{GD} in the context of ``formal differential geometry". A general definition
of a Poisson structure on infinite jets was proposed by the author in \cite{Vham} but its 
extension to general diffieties appeared to be a not very trivial task. 

More precisely, while the necessary definition of multivectors in secondary calculus, 
sometimes also called variational multivectors, is a natural generalization of Definition 
\ref{def-sym}, some technical aspects of the related  Schouten bracket mechanism are
to be still elaborated. See \cite{V2,KeKrVe,GKV} for further results.

On the other hand, in the context of integrable systems numerous concrete Poisson 
structures were revealed. Among them the bi-hamiltonian ones deserves a special 
mention (see \cite{Mag}). So, arises the question of their classification.
In the finite-dimensional case the famous Darboux 
lemma tells that symplectic manifolds or, equivalently, nondegenerate Poisson structures
of the same dimension are locally equivalent. ``What is its analogue in field theory?" is a
good question, which, at first glance, seems to be out of place as many known examples 
show. Nevertheless, by substituting ``coverings" for ``diffeomorphisms" in the formulation 
of this lemma and observing that these two notions are locally identical for
finite-dimensional manifolds we get some satisfactory results. Namely, all Poisson 
structures explicitly described up to now on infinite jets are obtained from a few 
models by passing to suitable coverings. See \cite{AstV}) for more details.

\subsection{Nonlocal symmetries.}\label{nlc-sym}
The first idea about nonlocal symmetries its takes origin at a seemingly technical fact. 
It was observed that the PDEs forming the KdV 
hierarchy are obtained from the original KdV equation $\Ee=\{u_t=uu_x+u+u_{xxx}\}$ 
by applying the so-called {\it recursion operator}. This operator 
$$
\Rr=D_x^2+\frac{2}{3}u+\frac{1}{3}u_xD_x^{-1}
$$
is not defined rigorously. By applying it to generating functions of symmetries of 
$\Ee$ one gets new ones that may depend on $D_x^{-1}u=\int udx$. A due rigor to this 
formal trick can be given by passing to a 1-dimensional covering 
$\Ee_{\infty}^{\prime}\to\Ee_{\infty}$ by adding to standard coordinates on $\Ee_{\infty}$ 
a new one $w$ such that $D_xw=u$ and $D_tw=u_{xx}+\frac{1}{2}u^2$. In this setting 
the above symmetries of $\Ee$ depending on $\int dx$, i.e., {\it nonlocal} ones, become
symmetries of $\Ee_{\infty}^{\prime}$ in the sense of Definition \ref{def-sym}, i.e.,  {\it local}
ones. This and other similar arguments motivate the following definition.
\begin{defi}\label{non-def}
A \emph{nonlocal symmetry} $($finite or infinitesimal$)$ of an equation $\Ee$ is a local symmetry 
of a diffiety $\Oo$, which covers $\Ee_{\infty}$. If $\tau:\Oo\to\Ee_{\infty}$ is a covering, then 
symmetries of $\Oo$ are called  \emph{$\tau$--symmetries} of $\Ee$.
\end{defi}
Similarly are defined nonlocal quantities of any kind. For instance, Poisson
structures in field theory discussed in Subsection \ref{ex-cover} are nonlocal with respect to
the original PDE/diffiety.

Let $\tau_i:\Oo_i\to\Ee_{\infty}, \,i=1,2,$ be two coverings of $\Ee_{\infty}$. A remarkable fact, 
which is due to I.~S.~Krasil'shchik, is that the Lie bracket of a  $\tau_1$--symmetry and a 
$\tau_2$--symmetry can be defined as a $\tau$--symmetry for a suitable covering 
$\tau:\Oo\to\Ee_{\infty}$ together with coverings 
$\tau_i^{\prime}:\Oo\to\Oo_i, \,i=1,2,$ such that $\tau=\tau_i\circ\tau_i^{\prime}$. The covering
$\tau$ is not defined uniquely. Nevertheless, this non-uniqueness can be resolved by passing
to a common  covering for ``all parties in question". The Jacobi identity as well as other 
ingredients of Lie algebra theory can be settled in a similar manner (see \cite{KVN,KVin}).
So, nonlocal symmetries of a PDE $\Ee$ form this {\it strange} Lie algebra, and this fact in 
turn confirms the validity of Definition \ref{non-def}.

Thus this definition incorporates all theoretically or experimentally known candidates for
symmetries of a PDE. Moreover, it  brings us to a new challenging question: 
\begin{quote}
{\it Symmetries of which object are the elements \\ of the above ``strange" Lie algebra\,?} 
\end{quote}
Indeed, this algebra can be considered not only as the algebra of nonlocal symmetries of 
the equation $\Ee$ but also as the symmetry algebra of any  equation
that covers $\Ee$. In other words, the question: What are partial differential equations?
arises again in this new context. But before we shall take a necessary look at the related 
problem of construction of coverings.

\subsection{Finding of coverings.}
The problem of how to find coverings of a given equation is key from both practical and 
theoretical points of view. At present we are rather far from its complete solutions. So,
below we shall illustrate the situation by sketching a direct method, which works well for 
PDEs in two independent variables and also supplies us with an interesting experimental 
material.  

Let $\Ee\subset J^k(E,n), \,\Oo=(\M,\Dd)$ and $\tau:\Oo\to \Ee_{\infty}$ be a covering. 
A $\tau$--projectable vector field $X\in \Dd$ is of the form $X=\bar{X}+V$ where
$\bar{X}\in\Cc_{\Ee}$ and $V$ is $\tau$--vertical, i.e., tangent to the fibers of $\tau$.
Locally $\tau$ can be represented as the projection $U\times W\to\Ee_{\infty}$ 
with $W$ being a pro-finite manifold and $U$ a domain in $\Ee_{\infty}$. If $U$
is sufficiently  small, then the restrictions $\bar{D}_i,  i=1,\dots,n,$ of the total derivatives 
$D_i$'s to $U$ span the distribution $\Cc_{\Ee}|_{U}$. The vector fields 
$\widehat{D}_i\in \Dd$ that project onto the $\bar{D}_i$ span $\Dd|_{\tau^{-1}(U)}$ and
$\widehat{D}_i=\bar{D}_i+V_i$ where $V_i$ is $\tau$--vertical. The Frobenius property of
$\Dd$ is equivalent to 
\begin{equation}\label{com-cov}
0=[\widehat{D}_i,\widehat{D}_j]\Leftrightarrow [\bar{D}_i,V_j]-[\bar{D}_j,V_i]+[V_i,V_j]=0, 
\quad 1\leq i<j\leq n.
\end{equation}
By inverting this procedure we get a method to search for coverings of $\Cc_{\Ee}$. Namely,
take a pro-finite manifold $W$ with coordinates $w_1,w_2,\dots$ and vector fields 
$V_i=\sum_r a_r\p/\p w_s$ on $U\times W$ with indeterminate coefficients 
$a_r\in\ci(U\times W)$. Any choice of these coefficients satisfying relations (\ref{com-cov})
defines a Frobenius distribution $\Span\{\widehat{D}_1,\dots,\widehat{D}_n\}$, which covers 
$\Cc_{\Ee}$. So, by resolving equations (\ref{com-cov}) with respect to the $a_r$ we get local
coverings of $\Cc_{\Ee}$. Many exact solutions of these equations for concrete PDEs of
interest can be found for $n=2$ and they reveal a very interesting structure, which we illustrate
with the following example.  
\begin{ex}\label{KdV-cov}
For the KdV equation $\Ee=\{u_t=uu_x+u_{xxx}\}$ we may take $t,x,u,u_x,u_{xx},\dots$ for
coordinates on $\Cc_{\Ee}$. Then the vector fields
$$
D_x\df\bar{D}_1=\frac{\p}{\p x}+\sum_{s=0}^{\infty}u_{s+1}\frac{\p}{\p u_{s}}, \quad
D_t\df\bar{D}_2=\frac{\p}{\p t}+\sum_{s=0}^{\infty}D_x^s(u_{3}+uu_{1})\frac{\p}{\p u_{s}},
$$
with $u_{s}=u_{x\dots x}  \;(s$--times) span  $\Cc_{\Ee}$. Put $V_x=V_1, \,V_t=V_2$ and 
look for solutions of \emph{(\ref{com-cov})} assuming that  $a_r=a_r(u,u_1,u_2,w_1,w_2,\dots)$ 
for simplicity. The result is worth to be reported in details. We have
\begin{equation}\label{WE-vf}
\begin{array}{l}
V_x=u^2A+uB+C, \medskip\\
V_t=2uu_2A+u_2B-u_1^2A+u_1[B,C]+\frac{2}{3}u^3A+\smallskip\\
\qquad + \frac{1}{2}(B+[B,[C,B]])+u[C,[C,B]]+D
\end{array}
\end{equation}
with $A,B,C,D$ being some fields on $W$ such that
\begin{equation}\label{WE-com}
\begin{array}{l}
\,[A,B]=[A,C]=[C,D]=0, \quad [B,D]+[C,[C,[C,B]]]=0, \bigskip\\
\,[B,[B,[B,C]]]=0, \quad[A,D]+\frac{3}{2}[B,[C,[C,B]]]=0 .
\end{array}
\end{equation}
This results tells that if we consider the Lie algebra  generated by four elements $A,B,C,D$, 
which are subject to the relations \emph{(\ref{WE-com})}, then any representation of this algebra by 
vector fields on a manifold $W$ gives a covering of $\Cc_{\Ee}$ associated with the vector fields 
\emph{(\ref{WE-vf})}.
\end{ex} 
\begin{rmk}
The Lie algebra 
defined by relations $($\ref{WE-com}$)$ ``mystically" 
appeared for the first time in the paper by H.~D.~Wahlquist and F.~B.~Estabrook \cite{WE}, 
in which they introduced the so-called prolongation structures. The fact that it is, as
 explained above, a necessary ingredient in the construction of coverings is due to the author.
\end{rmk} 
The reader will find many other examples of this kind together with related nonlocal symmetries, conservation laws, recursion operators, B\"acklund transformations, etc. in \cite{KVN,KVin,KK}.

\subsection{But what really are PDEs?} 
Now we can turn back to the question posed at the end of
Subsection \ref{nlc-sym}. Recall that the possibility to commute nonlocal symmetries of a PDE 
$\Ee$ living in different coverings of $\tau_i:\Oo_i\to \Ee_{\infty}, \,i=1,\dots,m,$ is ensured by 
the existence of a common covering diffiety $\Oo$, i.e., a system of coverings 
$\tau_i^{\prime}:\Oo\to \Oo_i$ such that $\tau=\tau_i\circ\tau_i^{\prime}$. So, in order 
to include into consideration all nonlocal symmetries we must consider ``all" coverings  
$\tau_{\a}:\Oo_{\a}\to\Ee_{\infty}$ of $\Ee_{\infty}$ as well as coverings $\Oo_{\a}\to\Oo_{\b}$. 
In this way we come to the category $\mathbf{Cobweb}\,\Ee$ of coverings of $\Ee_{\infty}$.
Then it is natural to call the {\it universal covering} of $\Ee$ the {\it terminal object} of  
$\mathbf{Cobweb}\,\Ee$. Denote this hypothetical universal covering by
$\tau_{\Ee}:\Oo_{\Ee}\to\Ee_{\infty}, \,\Oo_{\Ee}=(\M_{\Ee},\Dd_{\Ee})$. 
Now it is easy to see that $\mathbf{Cobweb}\,\Ee=\mathbf{Cobweb}\,\Ee^{\prime}$ if and 
only if there is a
common covering diffiety $\Oo$,  \,$\Ee_{\infty}\lw\Oo\rw\Ee^{\prime}_{\infty}$.
In other words, $\Ee$ and $\Ee^{\prime}$ are related by a B\"acklund transformation (see
Subsection \ref{ex-cover}). Recalling that coverings present inversions of differential 
operators we can trace the following analogy with algebraic geometry:  \smallskip
\begin{itemize}
\item []  {\it Affine algebraic variety associated with an algebraic equation} 
$\;\Rightarrow \;\Ee_{\infty}$. \smallskip
\item []  {\it Birational transformations connecting two affine varieties} $\;\Rightarrow$ 
\it {B\"acklund transformations} .\smallskip
\item [] {\it The field of rational functions on an affine variety} $\;\Rightarrow \;\Oo_{\Ee}$.
\end{itemize} \smallskip
This analogy becomes a tautology if one considers algebraic varieties as PDEs in {\it zero
independent variables}. Indeed, any DO in this case is of zero order, i.e., multiplications by
a function, and hence the inversion of such a DO is  the division by this function.

Unfortunately, the universal covering understood as a terminal object of a category is not 
sufficiently constructive to work with. 
However, we have some indication of how to proceed. From the theoretical side, the indication
is to look for an analogue of the fundamental group in the category of diffieties in order to 
construct the universal covering. By taking into account that we deal with infinitesimal
symmetries it would be more adequate to look for the {\t infinitesimal fundamental group}, 
i.e., for the {\it fundamental Lie algebra} of the diffiety $\Ee_{\infty}$. On the other side, 
this idea is on an ``experimental" ground. Namely, the Lie algebra associated with a 
Wahlquist-Estabrook prolongation structure (see Example \ref{KdV-cov}) is naturally 
interpreted as the universal algebra for a special class of coverings. 

In this connection a very interesting result by S.~Igonin should be mentioned. In \cite{Igo}
he constructed an object which possesses basic properties of the fundamental algebra
for a class of PDEs in two independent variables. Moreover, on this basis he succeeded 
to prove the non-existence of B\"acklund transformations connecting some integrable PDE's, 
for instance, the KdV equation and the Krichever-Novikiov equation.

Thus the question: What are PDE's ? continues to resist well, and the reader may see that
this is a highly nontrivial conceptual problem. Yet though universal coverings of diffieties
(if they exist\,!) point at a plausible answer, a good bulk of work should be done in order 
to put these ideas on a firm ground.   

\section{A couple of words about secondary calculus.}
In these pages we, first, tried to attract attention to two intimately related questions: ``what 
are symmetries of an object?" and ``what is the object itself?". They form something like an
electro-magnetic wave when one of them induces the other and vice versa. Probably, this
dynamical form is the most adequate adaptation of the background ideas of the Erlangen 
program to realities of present-day mathematics. The launch of such a wave in the area of
nonlinear partial differential equations was the inestimable contribution of S.~Lie to 
modern mathematics as it is now clearly seen in the hundred-years retrospective. 

In the above picture of the post-Lie phase of propagation of this wave we did not touch 
such fundamental questions as what are general tensor fields, connections, differential 
operators, etc, on the ``space of all solutions" of a given PDE, i.e., on the corresponding
diffiety. They all together form what we call {\it secondary calculus}. It turns out that any
natural notion or construction of the standard ``differential mathematics" has an analogue
in secondary calculus, which is referred to by adding the adjective ``secondary". In these
terms (higher) symmetries of a PDE $\Ee$ are nothing but {\it secondary vector fields} on 
$\Ee_{\infty}$. Surprisingly, all secondary notions are cohomology classes of suitable
natural complexes of differential operators, one of which, the jet-Spencer complex, was
discussed in Section \ref{inter-dif}. For the whole picture see \cite{V2}.

To illustrate this point we shall give some details on secondary differential forms. They
constitute the first term of the {\it $\Cc$--spectral sequence}, which is defined as follows.
Let $\Oo=(\M,\Dd)$ be a diffiety and 
$\Dd\Lambda(\Oo)=\oplus_{i\geq 0}\Dd\Lambda^i(\Oo)$ the ideal of differential forms 
on $\M$ vanishing on the distribution $\Dd$. This ideal is differentially closed and its
powers $\Dd^k\Lambda(\Oo)$ form a decreasing filtration of $\Lambda(\Oo)$. The
$\Cc$--spectral sequence $\{E^{p,q}_r(\Oo),d^{p,q}{\Oo}\}$ is the spectral sequence 
associated with this filtration. 
By definition, the space of secondary differential forms of degree $p$ is the graded object
$\oplus_{q=0}^{n}E^{p,q}_1(\Oo)$ and $d_1$ is the secondary exterior differential. 
Note that a smooth fiber bundle may be naturally viewed as a diffiety and the corresponding
$\Cc$--spectral sequence is identical to the Leray-Serre spectral sequence of  this bundle.  

Nontrivial terms of the $\Cc$--spectral sequence are all in the strip $0\leq q\leq n, \,p\geq 0$
with $n=\dim \Dd$, and $E_1^{0,q}(\Oo)=\bar{H}^q(\Oo)$ (horizontal de Rham cohomology
of $\Oo$, see Subsection \ref{hor}). Below we write simply $E_r^{p,q}$ for $E^{p,q}_r(\Oo)$ 
if the context does not allow a confusion. Also recall that $\Cc$--differential DOs are those 
that admit restrictions to integral submanifolds of $\Dd$.

The following proposition illustrates the fact that the calculus of variations is just an element 
of the calculus of secondary differential forms.
\begin{proc}\label{C-jet}
Let $\Oo=J^{\infty}(E,n)$. Then
\begin{enumerate}
\item If $E_1^{p,q}$ is nontrivial, then either $p=0$ or $q=n$  \emph{(``one line theorem")}.
\smallskip
\item $E_1^{0,q}=H^q(J^1(E,n))$,  if $q<n$, and $E_1^{0,n}$ is composed of variational
functionals $\int\omega\bar{d}x_1\wedge\dots\wedge\bar{d}x_n$.\smallskip
\item $d^{0,n}_1$ is the \emph{Euler operator} of the calculus of variations:
$$ E_1^{0,n}=\bar{H}^n(J^{\infty}(E,n))\ni
\int\omega\bar{d}x_1\wedge\dots\wedge\bar{d}x_n \;
\stackrel{d^{0,n}_1}{\longmapsto} \;\ell_{\omega}^*(1)\in\widehat{\vk}
$$
where $\widehat{\vk}\df\Hom_{\Ff}(\vk,\bar{\Lambda}^n(J^{\infty}(E,n)))$ and $\ell_{\omega}^*$
stands for the adjoint to $\ell_{\omega}$ $\Cc$--differential operator.\smallskip
\item $E_1^{2,n}=\Cc\Diff^{\mathrm{alt}}(\vk,\widehat{\vk})=\{$skew-self-adjoint $\Cc$--differential
operators from $\vk$ to $\widehat{\vk}\}$,  and 
$$
d^{1,n}_1: \widehat{\vk}\ni \Psi \;\longmapsto
\;\ell_{\Psi}^*-\ell_{\Psi}^*\in\Cc\Diff^{\mathrm{alt}}(\vk,\widehat{\vk}) .
$$
\item $E_2^{p,n}=H^{p+n}(J^1(E,n))$ and, in particular, the complex 
$\{E_1^{p,n},d_1^{p,n}\}_{p\geq 0}$ is locally acyclic.
\end{enumerate} 
\end{proc}
\noindent The reader will find a similar description of the terms $E_1^{p,n}$ and the 
differentials $d_1^{p,n}$ for $p>2$ in \cite{V-C-spec,V2}.

If $\Oo=\Ee_{\infty}$, then the terms $E_1^{0,q}(\Oo)$ present various conserved quantities of
the equation $\Ee$. For instance, the Gauss electricity conservation law is an element of
$E_1^{0,2}(\Ee_{\infty})$ for the system of Maxwell equations $\Ee$. The term 
$E_1^{0,n-1}(\Ee_{\infty})$ is composed of standard  {\it conservation laws} of a PDE $\Ee$,
which are associated with {\it conserved densities}. In this connection we have
\begin{proc}\label{2-lines}
Let $\Ee$ be a determined system of PDEs and $\Cc\Ll(\Ee)\df E_1^{0,n-1}(\Ee_{\infty})$ 
the vector space of conservation laws for $\Ee$. Then
\begin{enumerate}
\item If $E_1^{p,q}$ is nontrivial, then either $p=0$ or $q=n-1,n$ \emph{(``two lines theorem")}.\smallskip
\item $\ker d_1^{0,n-1}=H^{n-1}(\Ee)$ \emph{(trivial conservation laws)}.
\item $E_1^{1,n-1}=\ker\ell_{\Ee}^*$ and $E_1^{1,n}=\mathrm{coker} \,\ell_{\Ee}^*$.
\end{enumerate} 
\end{proc}
$\Upsilon=d_1^{0,n-1}(\Omega)$ is called the {\it  generating function} of a conservation law 
$\Omega\in \Cc\Ll(\Ee)$. Assertion (2) of  Proposition  \ref{2-lines} tells that a conservation 
law is uniquely defined by its generating function up to a trivial one. Moreover, by assertion 
(3) of this proposition, generating functions are solutions of the equation
$\ell_{\Ee}^*\Upsilon=0$, and this is the most efficient known method for finding 
conservation laws (see \cite{V3,KVin,KVe}).  

Propositions  \ref{C-jet} and \ref{2-lines} unveil the nature of the classical Noether theorem.
Namely, by assertions (3) and (4) of Proposition  \ref{C-jet}, the Euler-Lagrange equation $\Ee$
corresponding to the Lagrangian   $\int\omega\bar{d}x_1\wedge\dots\wedge\bar{d}x_n$  is
$\Psi=0$ with $\Psi=\ell_{\omega}^*(1)$ and $\ell_{\Ee}=\ell_{\Ee}^*$. In other words,    
Euler-Lagrange equations are {\it self-adjoint}. So, in this case the equation $\ell_{\Ee}^*=0$ 
whose solutions are generating functions of conservation laws of $\Ee$ (assertion (3) of Proposition  
\ref{2-lines}) coincides with the equations  $\ell_{\Ee}=0$ whose solutions are generating functions 
of symmetries of $\Ee$  (formula (\ref{eq.sym})). Moreover, we see that this relation between
symmetries and conservation laws takes place for a much larger than the Euler-Lagrange 
class of PDEs, namely, the class of {\it conformally self-adjoint} equations : 
$\ell_{\Ee}^*=\lambda\ell_{\Ee}, \,\lambda\in\Ff_{\Ee}$.

All natural relations between vector fields and differential forms such as Lie derivatives, insertion
operators, etc survive at the level of secondary calculus in the form of some relations between
the horizontal jet-Spencer cohomology and the first term of the $\Cc$--spectral sequence. Also, 
a morphism of diffieties induces a pull-back homomorphism of $\Cc$--spectral sequences. In
particular, this allows to define {\it nonlocal conservation laws} of a PDE ${\Ee}$ as conservation 
laws of diffieties that cover ${\Ee}_{\infty}$. These are just a few of numerous facts that show
high self-consistence of secondary calculus and its adequacy for needs of physics and 
mechanics.
\begin{rmk}
The $\Cc$--spectral sequence was introduced by the author in \cite{V-C-D}. It was preceded
by some works by various authors on the inverse problem of calculus of variations and the 
resolvent of the Euler operator $($or the Lagrange complex$)$. These works may now be 
seen as results about the $\Cc$--spectral sequence for $\Oo=J^{\infty}(\pi)$ $($see, for 
instance, \cite{Kuper, Tul}$)$. If $\Ee\subset J^{k}(\pi)$, then the first term of the $\Cc$--spectral 
sequence for $\Ee_{\infty}\subset J^{\infty}(\pi)$ acquires the second differential coming 
from the spectral sequence of the fiber bundle $\pi_{\infty}:J^{\infty}(\pi)\rw M$ and it becomes
the \emph{variational bicomplex} associated with $\Ee$. This local interpretation of the 
$\Cc$--spectral sequence is due to T.~Tsujishita \cite{Tsu}, who described these two 
differentials in a semi-coordinate manner.
\end{rmk}

\section{New language and new barriers.} 
In the preceding pages we were trying to show that a pithy general theory of PDEs exists
and to give an idea about the new mathematics that comes into light when developing
this theory in 
a systematic way. Even now this young theory provides many new instruments allowing 
to discover new features and facts about well-known and for long time studied PDE's in 
geometry, mechanics and  mathematical physics. The theory of singularities of solutions 
of PDs sketched in Section \ref{high-cont-str} is an example of this to say nothing about 
symmetries, conservation laws, hamiltonian structures and other more traditional aspects. 
Moreover, numerous possibilities, which are within one arm's reach, are still waiting to be duly
elaborated simply because of a lack of workmen in this new area. This situation is to a great 
extent due to a {\it language barrier}, since
\begin{quote} 
{\it the specificity of the general theory of PDEs is that it cannot be systematically 
developed in all its aspects on the basis of the traditionally understood differential calculus}. 
\end{quote}
Indeed, one very soon loses the way 
by performing exclusively direct manipulations with coordinate-wise descriptive definitions of 
objects of differential calculus, especially if working on such infinite-dimensional objects as 
diffieties. By their nature, these descriptive definitions cannot be applied to various situations 
when some kind of singularities or other nonstandard situations occur naturally. Not less 
important is that descriptive definitions give no idea about natural relations between objects 
of differential calculus. Typical questions that can in no way be neglected when dealing with 
foundations of the theory of PDE's are: ``What are tensor fields on manifolds with singularities, 
or on pro-finite manifolds, or what are tensor fields respecting a specific structure on 
a smooth manifolds", etc. This kind of questions becomes much more delicate when working 
with diffieties.  

All these questions can be answered by analyzing  why and how the traditional 
differential calculus of Newton and Leibniz became a natural language of classical physics
(including geometry and mechanics). Since the fundamental paradigm of classical physics 
states that existence means observability and vice versa, the first step in this analysis must 
be a due mathematical formalization of the observability mechanism in classical physics. 

We do that by assuming that from a mathematical point of view a classical physical laboratory 
is the {\it unitary algebra} $A$ over $\dR$ {\it generated by measurement instruments} installed 
in this laboratory and called the {\it algebra of observables}. A state of an observed object  
is interpreted as a homomorphism $h:A\to \dR$ of $\dR$--algebras ($\equiv$ ``readings of all 
instruments"). Hence the variety of all states of the system is identified with the real spectrum 
$\Spec_{\dR}A$ of $A$. 
The validity of this formalization of the classical observation mechanism is confirmed by the 
fact that all aspects of classical physics are naturally and, even  more, elegantly expressed in 
terms of this language. Say, one of the simplest necessary concepts, namely, that of velocity 
of an object at a state $h\in\Spec_{\dR}A$ is defined as a tangent vector to $\Spec_{\dR}A$ 
at the ``point" $h$, i.e., 
as an $\dR$--linear map $\xi:A\to\dR$ such that $\xi(ab)=h(a)\xi(b)+h(b)\xi(a)$. So, velocity 
is a particular first order DO over the algebra $A$ of observables in the sense of Definition
\ref{def-DO}. In this case $P=A$ and $Q=\dR$ as an $\dR$--vector space with the 
$\dR$--module product $a\star r\df h(a)r, \,a\in A, \,r\in\dR$. The reader will find other simple examples of this kind in an elementary introduction to the subject \cite{Nestr}.

Thus, by formalizing the concept of a classical physical laboratory as a commutative algebra,
we rediscover differential calculus in a new and much more general form. The next question is :
``What is the structure of this new language and what are its informative capacities?"  In the 
standard approach the zoo of various structures and constructions in modern differential and algebraic geometry, mechanics, field theory, etc. that are based on differential calculus seems 
not to manifest any regularity. Moreover, numerous questions like ``why do skew symmetric 
covariant tensors, i.e., differential forms, possess a natural differential $d$, while the symmetric 
ones do not" cannot be answered within this approach. On the contrary, in the framework 
of differential calculus over commutative algebras all these ``experimental materials" are nicely 
organized within a scheme composed of {\it functors of differential calculus connected by natural transformations} and {\it the objects that represent them} in various categories of modules 
over the ground algebra. 

The reader will find in a series of notes \cite{VV} various examples illustrating what  one can discover by analyzing the question ``what is the conceptual definition of covariant tensors". From the last three
notes of this series he can also  get an idea on the complexity of the theory of iterated differential 
forms and, in particular, tensors, in secondary calculus.

It should be especially mentioned that new views, instruments and facts coming from the general theory of PDEs and related mathematics offer not only new perspectives for many branches of contemporary mathematics and physics but at the same time put in question some popular 
current approaches and expectations ranging from algebraic geometry to QFT. Unfortunately,
there is too much to say in order to present the necessary reasons in a satisfactory manner.\newline

We conclude by stressing that
\begin{quote}
{\it The complexity and the dimension of problems in general theory of PDEs are so high that 
a new organization of mathematical research similar to that in experimental physics is absolutely indispensable.}
\end{quote}
Unfortunately, the dominating mentality and the ``social organization" of the modern 
mathematical community seems not to be sufficiently adequate to face this challenge.

{\bf Acknowledgements.} The author is very grateful to A.~B.~Sossinsky who 
moved him to write a ``popular" paper dedicated to modern general theory 
of PDEs and addressed to a general audience, and to the Editors for their 
kind invitation to include such a paper in the present volume. Additional thanks 
to Katya Vinogradova for correcting author's ``continental English" and to Athanase
Papadopoulos for his improving the text suggestions.

\end{document}